\documentclass[a4paper,12pt,reqno]{amsart}
\usepackage{amsfonts}
\usepackage{amsmath}
\usepackage{amssymb}
\usepackage[a4paper]{geometry}
\usepackage{mathrsfs}
\usepackage{csquotes}
\usepackage{enumitem}
\usepackage{dirtytalk}
\usepackage{comment}
\usepackage{nccmath}

\usepackage[colorlinks]{hyperref}
\renewcommand\eqref[1]{(\ref{#1})} 
\numberwithin{equation}{section}
\theoremstyle{plain}
\newtheorem{thm}{Theorem}[section]

\newtheorem{cor}[thm]{Corollary}
\newtheorem{lem}[thm]{Lemma}

\theoremstyle{definition}

\newtheorem{rem}[thm]{Remark}

\newcommand\norm[1]{\left\lVert#1\right\rVert}

\newcommand*{\QEDA}{\null\nobreak\hfill\ensuremath{\square}}

\def\e[#1]{{\textrm{e}}^{#1}}

\begin{document}

   \title[General weighted Hardy type inequalities]
 {Refined general weighted $L^p$-Hardy and Caffarelli-Kohn-Nirenberg type inequalities and identities related to the Baouendi-Grushin operator}

\author[N. Yessirkegenov]{Nurgissa Yessirkegenov}
\address{
  Nurgissa Yessirkegenov:
  \endgraf
  KIMEP University, Almaty, Kazakhstan
     \endgraf
  {\it E-mail address} {\rm nurgissa.yessirkegenov@gmail.com}
  }

  \author[A. Zhangirbayev]{Amir Zhangirbayev}
\address{
  Amir Zhangirbayev:
 \endgraf
   SDU University, Kaskelen, Kazakhstan
  \endgraf
  and 
  \endgraf
  Institute of Mathematics and Mathematical Modeling, Kazakhstan
   \endgraf
  {\it E-mail address} {\rm
amir.zhangirbayev@gmail.com}
  }

\thanks{This research is funded by the Committee of Science of the Ministry of Science and Higher Education of the Republic of Kazakhstan (Grant No. AP23490970).}


     \begin{abstract} In this paper, we present a sufficient condition on a pair of nonnegative weights $v$ and $w$ such that we have a general weighted $L^{p}$-Hardy type identity. The result, for a certain choice of weights, gives weighted $L^{p}$-Hardy type inequalities and identities with explicit remainder terms, thereby improving previously known results. Furthermore, we obtain the corresponding general weighted Caffarelli-Kohn-Nirenberg type inequality with remainder terms, which, as a result, imply Heisenberg-Pauli-Weyl type inequalities.

     \end{abstract}
     \maketitle

\section{Introduction}

First, let us recall the classical $L^p$-Hardy inequality on the usual Euclidean space $\mathbb{R}^n$: let $1<p<n$ and $f\in C_0^{\infty}(\mathbb{R}^n\backslash\{0\})$. Then, we have
\begin{align}\label{class_multdim_hardy}
\int_{\mathbb{R}^n}|\nabla f|^pdz\geq\left(\frac{n-p}{p}\right)^p\int_{\mathbb{R}^n}\frac{|f|^p}{|z|^p}dz, 
\end{align}
where $\left(\frac{n-p}{p}\right)^p$ is sharp but not achieved. The literature on the classical Hardy inequality (\ref{class_multdim_hardy}) is enormous, and we are unable to list all the related articles in this paper. Therefore, we refer only to standard monographs such as \cite{KMP06, KMP07, BEL15}. Since in this paper we are mainly interested in the results regarding the Baouendi-Grushin operator, from this point onward, let us recall results only in this direction.

The first extension of the sharp Hardy inequality (\ref{class_multdim_hardy}, for $p=2$, to the Baouendi-Grushin vector fields) was done by Garofalo \cite{Gar93}. To be more precise, the following inequality was obtained for all $f \in C_0^{\infty}\left(\mathbb{R}^m \times \mathbb{R}^k \backslash\{(0,0)\}\right)$:  
\begin{align}\label{Gar}
\int_{\mathbb{R}^n}|\nabla_{\gamma}f|^2 dz\geq \left(\frac{Q-2}{2}\right)^2 \int_{\mathbb{R}^n}\frac{|x|^{2 \gamma}}{\rho^{2\gamma+2}}|f|^2 dz
\end{align}
with $n=m+k$, $Q=m+(1+\gamma)k$ and $\rho=\left(|x|^{2\gamma+2}+(1+\gamma)^2|y|^2\right)^{\frac{1}{2\gamma+2}}$. In this context, $\nabla_{\gamma}=(\nabla_{x},|x|^{\gamma}\nabla_{y})$, where $\nabla_{x}$ is the gradient operator in the variable $x$, and $\nabla_{y}$ is the gradient in the variable $y$. The inequality (\ref{Gar}) recovers the $L^{2}$ version of (\ref{class_multdim_hardy}) for $\gamma=0$.

Following Garofalo's influential work \cite{Gar93}, there has been a considerable amount of research focused on developing inequalities of Hardy-type related to the Baouendi-Grushin operator, as demonstrated in publications such as \cite{D'A04, NCH04, D'A05, DGN10, SY12, SJ12, Kom15, yang2015improved, KY18, LRY19, suragan2023sharp, d2024weighted, GJR24}. 

For example, in \cite[Formula (3.1) of Theorem 3.1]{D'A04}, D'Ambrosio was the first to obtain weighted $L^p$-versions of (\ref{Gar}): let $1<p<\infty$ and $\Omega$ be an open subset of $\mathbb{R}^n$. Then, we have
\begin{align}\label{D'A04}
\int_{\Omega}|x|^{\beta-\gamma p} \rho^{(1+\gamma) p-\alpha} \left|\nabla_\gamma f\right|^p dz \geq c^{p}_{Q,p,\alpha,\beta} \int_{\Omega} \frac{|x|^\beta}{\rho^\alpha} |f|^p dz
\end{align}
for all $f\in D^{1,p}_{\gamma}(\Omega,|x|^{\beta-\gamma p}\rho^{(1+\gamma)p-\alpha})$, $p>1$, $k,m\geq1$, $\alpha,\beta\in\mathbb{R}$ such that $m+(1+\gamma)k>\alpha-\beta-p$ and $m>\gamma p-\beta$. The space $D^{1,p}_{\gamma}(\Omega, \omega)$ denotes the closure of $C^{\infty}_{0}(\Omega)$ in the norm $(\int_{\Omega}^{}|\nabla_{\gamma}f|^{p}\omega dz)^{\frac{1}{p}}$ for $\omega\in L^{1}_{loc}(\Omega)$ with $\omega>0$ a.e. on $\Omega$. The constant $c^{p}_{Q,p,\alpha,\beta}$ in (\ref{D'A04}) is equal to $\left(\frac{Q+\beta-\alpha}{p}\right)^p$ and sharp when $0\in\Omega$. 

Around the same time, Niu, Chen and Han \cite{NCH04} utilized the Picone identity to prove a variety of Hardy-type inequalities for the whole space, the open pseudo-ball and the external domain of the pseudo-ball. 

Later, D'Ambrosio \cite[Formula (3.8) of Theorem 3.1]{D'A05} showed Hardy-type inequalities related to quasilinear second-order degenerate elliptic differential operators. Specifically, the following inequality was obtained: let $\alpha<-1$, $p=Q>1$ and $R>0$. Then, we have
\begin{align}\label{D'A05}
\int_{B_{R}}\left(\log\frac{R}{\rho}\right)^{\alpha+p}|\nabla_{\gamma}f|^{p}dz\geq\left(\frac{|\alpha+1|}{p}\right)^{p}\int_{B_{R}}\left(\log \frac{R}{\rho}\right)^{\alpha}\frac{|x|^{\gamma p}}{\rho^{\gamma p+p}}|f|^{p}dz
\end{align}
for all $f\in D^{1,p}_{\gamma}\left(B_{R},\left(\log\frac{R}{\rho}\right)^{\alpha+p}\right)$ and $B_{R}=\{z\in\mathbb{R}^n:\rho(z)<R\}$. Furthermore, the constant $\left(\frac{|\alpha+1|}{p}\right)^p$, in (\ref{D'A05}), is sharp.

Shen and Jin, in \cite{SJ12}, obtained Hardy-Rellich type inequalities with the optimal constant by a direct and simple approach. Then, in \cite{Kom15}, Kombe proved weighted Hardy-type inequalities with remainder terms, which play a key role in establishing improved Rellich type inequalities. Yang, Su and Kong, in \cite{yang2015improved}, derived improved Hardy inequalities on bounded domains containing the origin. In \cite{LRY19}, Laptev, Ruzhansky and the first author of this paper obtained a refinement of the Hardy inequalities and derived weighted Hardy type inequalities for the quadratic form of the magnetic Baouendi-Grushin operator for the magnetic field of Aharonov-Bohm type. Also, Suragan and the first author of this paper, in \cite{suragan2023sharp}, established a sharp remainder formula for the Poincaré inequality with a straightforward proof without using the variational principle. Recently, D'Arca in \cite{d2024weighted} derived Poincaré inequality and Hardy improvements related to some degenerate elliptic differential operators, including the Baouendi-Grushin operator. Then, Ganguly, Jotsaroop and Roychowdhury \cite{GJR24} showed Hardy, Hardy-Rellich, and Rellich identities and inequalities with sharp constants via spherical vector fields for $\gamma=1$ and $k=1$ in the setting of $L^{2}$. On top of that, such inequalities have also been investigated for other sub-elliptic operators of various types (see, e.g. \cite{GL90, Gar93, D'A04, D'A05, kombe2006hardy, DGN10, d2024weighted}).

Nevertheless, the result that we are most interested in is the one obtained by Kombe and Yener \cite{KY18}. There they derived a sufficient condition on a pair of non-negative weight functions $v$ and $w$ in $\mathbb{R}^{n}$ so that the general weighted Hardy type inequality with a remainder term is satisfied: let $v \in C^1\left(\mathbb{R}^n\right)$ and $w \in L_{\text {loc }}^1\left(\mathbb{R}^n\right)$ be nonnegative functions and $\phi \in$ $C^{\infty}\left(\mathbb{R}^n\right)$ be a positive function satisfying the differential inequality
\begin{align*}
-\nabla_\gamma \cdot\left(v\left|\nabla_\gamma \phi\right|^{p-2} \nabla_\gamma\phi\right) \geq w\phi^{p-1}
\end{align*}
a.e. in $\mathbb{R}^n$. There is a positive number $c_p=c(p)$ such that; if $p \geq 2$, then
\begin{align}\label{kom1}
\int_{\mathbb{R}^n} v\left|\nabla_\gamma f\right|^p dz \geq \int_{\mathbb{R}^n} w|f|^p dz+c_p \int_{\mathbb{R}^n} v\left|\nabla_\gamma \frac{f}{\phi}\right|^p \phi^p dz
\end{align}
and if $1<p<2$, then
\begin{align}\label{kom2}
\int_{\mathbb{R}^n} v\left|\nabla_\gamma f\right|^p dz \geq \int_{\mathbb{R}^n} w|f|^p dz+c_p \int_{\mathbb{R}^n} \frac{v\left|\nabla_\gamma \frac{f}{\phi}\right|^2 \phi^2}{\left(\left|\frac{f}{\phi} \nabla_\gamma \phi\right|+\left|\nabla_\gamma \frac{f}{\phi}\right| \phi\right)^{2-p}} dz
\end{align}
for all real-valued functions $f \in C_0^{\infty}\left(\mathbb{R}^n\right)$.

The main purpose of this paper is to provide similar results, but for complex-valued functions and nonnegative explicit remainder terms. In particular, we obtain the following refined general weighted $L^{p}$-Hardy type identity: let $1<p<\infty$ and let $\Omega\subseteq\mathbb{R}^{m+k}$ be an open set such that the integrals below make sense. Let $v\in C^{1}(\Omega\backslash\Sigma)$ and $w\in L^{1}_{loc}(\Omega\backslash\Sigma)$ be nonnegative functions satisfying the condition (in the weak sense)
\begin{align}\label{main condt intro}
0\leq\phi:=\nabla_{\gamma}\cdot \left(w^{\frac{p-1}{p}}v^{\frac{1}{p}}\frac{\nabla_{\gamma}\rho}{|\nabla_{\gamma}\rho|}\right)-pw
\end{align}
a.e. in $\Omega\backslash\Sigma$. Then, for all complex-valued $f\in C^{\infty}_0(\Omega\backslash\Sigma)$, we have
\begin{align}\label{main res intro}
\int_{\Omega}^{}v\frac{|\nabla_{\gamma}\rho\cdot\nabla_{\gamma}f|^{p}}{|\nabla_{\gamma}\rho|^{p}}dz=\int_{\Omega}^{}w|f|^{p}dz+\int_{\Omega}^{}C_{p}(\xi,\eta)dz+\int_{\Omega}\phi|f|^{p}dz,
\end{align}
where the functional $C_p(\cdot,\cdot)$ is given by 
\begin{align}\label{cp formula intro}
C_p(\xi,\eta):=|\xi|^p-|\xi-\eta|^p-p|\xi-\eta|^{p-2}\textnormal{Re}(\xi-\eta)\cdot\overline{\eta}\geq0
\end{align}
and
\begin{align}\label{main notation intro}
\xi:=v^{\frac{1}{p}}\frac{\nabla_{\gamma}\rho \cdot \nabla_{\gamma}f}{|\nabla_{\gamma}\rho|},\quad \eta:=v^{\frac{1}{p}}\frac{\nabla_{\gamma}\rho \cdot \nabla_{\gamma}f}{|\nabla_{\gamma}\rho|} + w^{\frac{1}{p}}f
\end{align}
with $\Sigma$ being the set of singular points of $v$ and $w$. 

The result (\ref{main res intro}) improves the left-hand side of (\ref{kom1}) and (\ref{kom2}) by replacing $\nabla_{\gamma}f$ with $\frac{\nabla_{\gamma}\rho\cdot\nabla_{\gamma}f}{|\nabla_{\gamma}\rho|}$. Here, we also have nonnegative explicit remainder terms. Moreover, (\ref{main res intro}) holds true for all complex-valued functions, whereas (\ref{kom1}) and (\ref{kom2}) are valid only for real-valued functions.  

A wide range of weighted Hardy-type identities can be derived by selecting appropriate weights $w$ and $v$ that satisfy the required condition (for more details, see Section \ref{subsection1}). For example, in the special case, when $v=|x|^{\beta-\gamma p}\rho^{p(1+\gamma)-\alpha}$ and $w=\left(\frac{Q+\beta-\alpha}{p}\right)^{p}\frac{|x|^{\beta}}{\rho^{\alpha}}$, in (\ref{main condt intro}), we have $\phi=0$, and a sharp remainder formula of (\ref{D'A04}) is obtained with $\frac{\nabla_{\gamma}\rho\cdot\nabla_{\gamma}f}{|\nabla_{\gamma}\rho|}$: let $1<p<\infty$, $m,k\geq1$ and $\alpha,\beta\in \mathbb{R}$ be such that $Q> \alpha-\beta$. Then, we have
\begin{multline}\label{dam ref}
\int_{\mathbb{R}^n}|x|^{\frac{\beta-\gamma p}{p}}\rho^{\frac{p(1+\gamma)-\alpha}{p}}\frac{|\nabla_{\gamma}\rho\cdot\nabla_{\gamma}f|^p}{|\nabla_{\gamma}\rho|^p}dz=\left(\frac{Q+\beta-\alpha}{p}\right)^{p}\int_{\mathbb{R}^n}\frac{|x|^{\beta}}{\rho^{\alpha}}|f|^pdz
\\+\int_{\mathbb{R}^n}C_{p}\biggl(|x|^{\frac{\beta-\gamma p}{p}}\rho^{\frac{p(1+\gamma)-\alpha}{p}}\frac{\nabla_{\gamma}\rho \cdot \nabla_{\gamma}f}{|\nabla_{\gamma}\rho|},|x|^{\frac{\beta-\gamma p}{p}}\rho^{\frac{p(1+\gamma)-\alpha}{p}}\frac{\nabla_{\gamma}\rho \cdot \nabla_{\gamma}f}{|\nabla_{\gamma}\rho|} \\+ \left(\frac{Q+\beta-\alpha}{p}\right)\frac{|x|^{\frac{\beta}{p}}}{\rho^{\frac{\alpha}{p}}}f\biggr)dz
\end{multline}
for all complex-valued $f\in C^{\infty}_0(\mathbb{R}^m\times\mathbb{R}^{k}\backslash\{(0,0)\})$. Here, we can use the results of Cazacu, Krej{\v{c}}i{\v{r}}{\'\i}k, Lam and Laptev \cite{cazacu2024hardy} to obtain the estimates of the remainder term from below for $p\geq2$ (see Corollary \ref{dam laptev}). Additionally, due to the recent results of Chen and Tang \cite{CT24} we can derive (\ref{dam ref}) with different sharp remainder terms for $1<p<2\leq n$ (see Corollary \ref{dam chinese}). Moreover, when $\beta=\gamma p$ and $\alpha=p(1+\gamma)$  in (\ref{dam ref}), the following identity is derived: let $1<p<\infty$ and $\gamma \geq0$. Then, for any complex-valued $f\in C^{\infty}_{0}(\mathbb{R}^m\times\mathbb{R}^{n}\backslash\{(0,0)\})$, there holds
\begin{multline}\label{gar ref intro lp}
\int_{\mathbb{R}^n}\frac{|\nabla_{\gamma}\rho\cdot\nabla_{\gamma}f|^{p}}{|\nabla_{\gamma}\rho|^{p}}dz=
\left(\frac{Q-p}{p}\right)^{p}\int_{\mathbb{R}^n}\frac{|x|^{\gamma p}}{\rho^{\gamma p + p}}|f|^{p}dz
\\+\int_{\mathbb{R}^n}C_{p}\left(\frac{\nabla_{\gamma}\rho \cdot \nabla_{\gamma}f}{|\nabla_{\gamma}\rho|},\frac{\nabla_{\gamma}\rho \cdot \nabla_{\gamma}f}{|\nabla_{\gamma}\rho|}+\left(\frac{Q-p}{p}\right)\frac{|x|^{\gamma}}{\rho^{\gamma+1}}f\right)dz.
\end{multline}
By setting $\gamma=0$ in (\ref{gar ref intro lp}), we obtain a refinement of the classical Hardy inequality (\ref{class_multdim_hardy}), which was also obtained, as a special case, by Kalaman and the first author of this article in \cite{kalaman2024cylindrical} under certain parameters: let $1<p<\infty$. Then, for any complex-valued functions $f\in C^{\infty}_{0}(\mathbb{R}^n\backslash\{0\})$, we have
\begin{multline*}
\int_{\mathbb{R}^n}\frac{|z\cdot\nabla f|^{p}}{|z|^{p}}dz=\left(\frac{n-p}{p}\right)^{p}\int_{\mathbb{R}^n}\frac{|f|^{p}}{|z|^{p}}dz\\+\int_{\mathbb{R}^n}C_{p}\left(\frac{z\cdot\nabla f}{|z|},\frac{z\cdot\nabla f}{|z|}+\left(\frac{n-p}{p}\right)f\right)dz.
\end{multline*} 
We refer to \cite{ruzhansky2017hardy} for the case of real-valued functions.

In addition, the identity (\ref{gar ref intro lp}) allows us to obtain a refinement of (\ref{Gar}) for $p=2$: let $\gamma\geq0$. Then, for any complex-valued $f\in C^{\infty}_{0}(\mathbb{R}^m\times\mathbb{R}^{n}\backslash\{(0,0)\})$, we have
\begin{multline*}
\int_{\mathbb{R}^n}\frac{|\nabla_{\gamma}\rho\cdot\nabla_{\gamma}f|^{2}}{|\nabla_{\gamma}\rho|^{2}}dz=
\left(\frac{Q-2}{2}\right)^{2}\int_{\mathbb{R}^n}\frac{|x|^{2\gamma}}{\rho^{2\gamma+2}}|f|^{2}dz
\\+\int_{\mathbb{R}^n}C_{2}\left(\frac{\nabla_{\gamma}\rho \cdot \nabla_{\gamma}f}{|\nabla_{\gamma}\rho|},\frac{\nabla_{\gamma}\rho \cdot \nabla_{\gamma}f}{|\nabla_{\gamma}\rho|}+\left(\frac{Q-2}{2}\right)\frac{|x|^{\gamma}}{\rho^{\gamma+1}}f\right)dz.
\end{multline*}

Another interesting application of our main result (\ref{main res intro}) is the refinement of the power logarithmic $L^p$-Hardy type inequality (\ref{D'A05}). That is, setting $v=\left(\log \frac{R}{\rho}\right)^{\alpha+p}$ and $w=\left(\frac{|\alpha+1|}{p}\right)^{p}\left(\log \frac{R}{\rho}\right)^{\alpha}\frac{|x|^{\gamma p}}{\rho^{\gamma p+p}}$, in (\ref{main condt intro}), gives the following sharp remainder formula of (\ref{D'A05}) with $\frac{\nabla_{\gamma}\rho\cdot\nabla_{\gamma}f}{|\nabla_{\gamma}\rho|}$: let $1<p<\infty$, $\alpha\in\mathbb{R}$ and $R>0$ be such that $\alpha+1<0$. Then, for all complex-valued  $f\in C^{\infty}_{0}(B_{R}\backslash\{(0,0)\})$, we have 
\begin{multline}\label{cor4 eq intro}
\int_{B_{R}}\left(\log\frac{R}{\rho}\right)^{\alpha+p}\frac{|\nabla_{\gamma}\rho\cdot\nabla_{\gamma}f|^p}{|\nabla_{\gamma}\rho|^p}dz=\left(\frac{|\alpha+1|}{p}\right)^{p}\int_{B_{R}}\left(\log \frac{R}{\rho}\right)^{\alpha}\frac{|x|^{\gamma p}}{\rho^{\gamma p+p}}|f|^{p}dz\\+\int_{B_{R}}C_{p}\Biggl(\left(\log\frac{R}{\rho}\right)^{\frac{\alpha+p}{p}}\frac{\nabla_{\gamma}\rho \cdot \nabla_{\gamma}f}{|\nabla_{\gamma}\rho|},\left(\log\frac{R}{\rho}\right)^{\frac{\alpha+p}{p}}\frac{\nabla_{\gamma}\rho \cdot \nabla_{\gamma}f}{|\nabla_{\gamma}\rho|} \\+\left(\frac{|\alpha+1|}{p}\right)\left(\log\frac{R}{\rho}\right)^{\frac{\alpha}{p}}\frac{|x|^{\gamma }}{\rho^{\gamma+1}}f\Biggr)dz
\\+\left(\frac{|\alpha+1|}{p}\right)^{p-1}(Q-p)\int_{B_{R}}\frac{|x|^{\gamma p}}{\rho^{\gamma p + p}}\left(\log \frac{R}{\rho}\right)^{\alpha+1}|f|^{p}dz.
\end{multline}
If we set $\gamma=0$, $\alpha=-n$ and $p=n$, in (\ref{cor4 eq intro}), then for $n\geq2$, we obtain the following scale-invariant critical version of the Hardy inequality derived by Ioku, Ishiwata and Ozawa \cite{ioku2015scale, ioku2016sharp}:
\begin{multline*}
\int_{B_{R}}\frac{|z\cdot\nabla f|^{n}}{|z|^{n}}dz=\left(\frac{n-1}{n}\right)^{n}\int_{B_R}\frac{|f|^{n}}{|z|^{n}\left(\log \frac{R}{|z|}\right)^n}dz\\+\int_{B_{R}}C_{n}\left(\frac{z\cdot\nabla f}{|z|},\frac{z\cdot\nabla f}{|z|}+\left(\frac{n-1}{n}\right)\frac{f}{|z|\left(\log \frac{R}{|z|}\right)}\right)dz
\end{multline*}
for all complex-valued $f\in C^{\infty}_{0}(B_{R}\backslash\{(0,0)\})$.

We also apply the main result to obtain the general weighted Caffarelli-Kohn-Nirenberg (CKN) type inequalities with remainder terms, which, in special cases, give Heisenberg-Pauli-Weyl (HPW) type inequalities. Before we state the results, let us recall the classical CKN inequality from \cite{caffarelli1984first}: let \( n \in \mathbb{N} \) and let $p,\ q,\ r,\ a,\ b,\ d,\ \delta \in \mathbb{R}$ such that \( p, q \geq 1,\ r > 0,\ 0 \leq \delta \leq 1 \), and
\begin{align*}
\frac{1}{p} + \frac{a}{n}, \quad \frac{1}{q} + \frac{b}{n}, \quad \frac{1}{r} + \frac{c}{n} > 0,
\end{align*}
where \( c = \delta d + (1 - \delta)b \). Then there exists a positive constant \( C \) such that
\begin{align}\label{ckn}
\left\| |x|^a |\nabla f| \right\|^{\delta}_{L^p(\mathbb{R}^n)} \left\| |x|^b f \right\|_{L^q(\mathbb{R}^n)}^{1 - \delta}\geq C\left\| |x|^c f \right\|_{L^r(\mathbb{R}^n)},
\end{align}
holds for all \( f \in C_0^\infty(\mathbb{R}^n) \) if and only if the following conditions hold:
\begin{align*}
\frac{1}{r} + \frac{c}{n} = \delta \left( \frac{1}{p} + \frac{a - 1}{n} \right) + (1 - \delta)\left( \frac{1}{q} + \frac{b}{n} \right),
\end{align*}
\begin{align*}
a - d \geq 0 \quad \text{if} \quad \delta > 0,
\end{align*}
\begin{align*}
a - d \leq 1 \quad \text{if} \quad \delta > 0 \quad \text{and} \quad \frac{1}{r} + \frac{c}{n} = \frac{1}{p} + \frac{a - 1}{n}.
\end{align*}
The classical CKN inequalities (\ref{ckn}) have been extended to a variety of settings, including the Baouendi-Grushin vector fields, general homogeneous Carnot groups and the Heisenberg group \cite{han2011class, ruzhansky2017horizontal, KY18, GJR24}. On the anisotropic analogues of the classical $L^{p}$-CKN inequalities, we refer to \cite{ozawa2019p, li2023anisotropic}. We also refer to other extended CKN inequalities as well as various versions and applications \cite{dolbeault2015rigidity, dolbeault2016rigidity, ruzhansky2017caffarelli, ruzhansky2017extended, lam2017sharp, ruzhansky2018extended, dong2018existence, cazacu2024caffarelli}.

Under specific choice of parameters, the classical CKN inequalities (\ref{ckn}) imply the classical HPW uncertainty principle \cite{weyl1950theory}:
\begin{align}\label{hpw}
\left( \int_{\mathbb{R}^n} |\nabla f|^2 \, dz \right)
\left( \int_{\mathbb{R}^n} |z|^2 |f|^2 \, dz \right)
\geq C \left( \int_{\mathbb{R}^n} |f|^2 \, dz \right)^2
\end{align}
for all \( f \in C_0^\infty(\mathbb{R}^n) \). The inequality (\ref{hpw}) mathematically expresses the fundamental physical principle that, in any quantum state, it is impossible to simultaneously determine both a particle’s position and momentum with complete accuracy.

Now we are ready to state our results regarding general weighted CKN type inequalities with remainder terms: let $1<p,q<\infty$, $0<r<\infty$ with $p+q\geq r$, $\delta \in [0,1] \cap \left[ \frac{r - q}{r}, \frac{p}{r} \right]$ and $b, c \in \mathbb{R}$. Assume that $\frac{\delta r}{p} + \frac{(1-\delta) r}{q} = 1$ and $c = \frac{\delta}{p} + b(1 - \delta)$. Then, for any complex-valued $f\in C^{\infty}_{0}(\mathbb{R}^m\times\mathbb{R}^{n}\backslash\{(0,0)\})$, we have
\begin{multline*}
\Biggl[\norm{v^{\frac{1}{p}}\frac{\nabla_{\gamma}\rho\cdot\nabla_{\gamma}f}{|\nabla_{\gamma}\rho|}}^{p}_{L^{p}(\mathbb{R}^n)}\\-\int_{\mathbb{R}^n}C_{p}\left(v^{\frac{1}{p}}\frac{\nabla_{\gamma}\rho \cdot \nabla_{\gamma}f}{|\nabla_{\gamma}\rho|},v^{\frac{1}{p}}\frac{\nabla_{\gamma}\rho \cdot \nabla_{\gamma}f}{|\nabla_{\gamma}\rho|} + w^{\frac{1}{p}}f\right)dz\Biggr]^{\frac{\delta}{p}}\norm{w^{b}f}^{1-\delta}_{L^{q}(\mathbb{R}^n)}\geq\norm{w^{c}f}_{L^{r}(\mathbb{R}^n)}.
\end{multline*}

After choosing appropriate weights and parameters, we obtain the following HPW type inequalities: let $\gamma\geq0$. Then, we have 
\begin{align}\label{almost hpw}
\left(\int_{\mathbb{R}^n}\frac{|\nabla_{\gamma}\rho\cdot\nabla_{\gamma}f|^{2}}{|\nabla_{\gamma}\rho|^{2}}dz\right)\left(\int_{\mathbb{R}^n}\frac{\rho^{2\gamma+2}}{|x|^{2\gamma}}|f|^{2}dz\right)\geq\left(\frac{Q-2}{2}\right)^{2}\left(\int_{\mathbb{R}^n}|f|^{2}dz\right)^{2}
\end{align}
for any complex-valued $f\in C^{\infty}_0(\mathbb{R}^m\times\mathbb{R}^{k}\backslash\{(0,0)\})$.

Substituting $\gamma=0$ and applying the Cauchy-Schwarz inequality, in (\ref{almost hpw}), we obtain the uncertainty principle (\ref{hpw}) with an explicit constant:
\begin{align*}
\left( \int_{\mathbb{R}^n} |\nabla f|^2 \, dz \right)
\left( \int_{\mathbb{R}^n} |z|^2 |f|^2 \, dz \right)
\geq \frac{(n-2)^2}{4} \left( \int_{\mathbb{R}^n} |f|^2 \, dz \right)^2.
\end{align*}

In Section \ref{prem}, we briefly review the essential definitions, notations and preliminary results. Section \ref{refined general} is devoted to establishing the main refined general weighted $L^p$‑Hardy type inequalities and identities for the range $1 < p < \infty$. Furthermore, we provide proofs of remainder estimates for $p\geq2$ and $1<p<2\leq n$ under certain conditions. Finally, in Section \ref{applications}, we explore applications of these results. Section \ref{subsection1} presents sharp refinements of various classical and weighted Hardy inequalities on $\mathbb{R}^n$, within $\rho$‑ball and involving logarithmic weights. In Section \ref{subsection2}, we derive general weighted CKN type inequalities with remainder terms, which in turn yield HPW type uncertainty principles.

\section{Preliminaries}\label{prem}

In this section, we will provide the notation and some preliminary results. Suppose $z=(x_{1},\ldots,x_{m},y_{1},\ldots,y_{k})$ or simply $z=(x,y)$ $\in\mathbb{R}^{m}\times\mathbb{R}^{k}$ with $m+k=n$ and $m,k\geq1$. The sub-elliptic gradient is the $n$-dimensional vector field given by
\[
\nabla_{\gamma}=(X_{1},\ldots,X_{m},Y_{1},\ldots,Y_{k}).
\]
Here, the components are defined as:
\[
X_{i}=\frac{\partial}{\partial{x_i}}, \quad i=1,\ldots,m, \quad Y_{j}=|x|^{\gamma}\frac{\partial}{\partial y_{j}}, \quad j=1,\ldots,k. 
\]
The Baouendi-Grushin operator on $\mathbb{R}^{n}$ is the operator
\[
\Delta_{\gamma} = \nabla_{\gamma}\cdot\nabla_{\gamma}=\Delta_{x}+|x|^{2\gamma}\Delta_{y},
\]
where $\Delta_{x}$ and $\Delta_{y}$ are Laplace operators in the variables $x\in\mathbb{R}^{m}$ and $y\in\mathbb{R}^{k}$, respectively. The Baouendi-Grushin operator is a sum of squares of $C^{\infty}$ vector fields satisfying the H\"ormander condition for even positive integers $\gamma$:
\begin{align}
\operatorname{rank} \operatorname{Lie}\left[X_1, \ldots, X_m, Y_1, \ldots, Y_k\right]=n \text {. } \nonumber
\end{align}
The anisotropic dilation attached to $\Delta_{\gamma}$ is defined as follows:
\begin{align}
\delta_\lambda(x, y)=\left(\lambda x, \lambda^{1+\gamma} y\right), \quad \lambda>0, \quad z=(x,y)\in\mathbb{R}^{m+k}. \nonumber 
\end{align}
The dilation's homogeneous dimension is given by
\begin{align}
Q=m+(1+\gamma) k.  \nonumber  
\end{align}
The Lebesgue measure's formula for alternating variables suggests that
\begin{align}
d \circ \delta_\lambda(x, y)=\lambda^Q d x d y. \nonumber
\end{align}
It is straightforward to check that
\begin{align}
X_i\left(\delta_\lambda\right)=\lambda \delta_\lambda\left(X_i\right), \quad Y_i\left(\delta_\lambda\right)=\lambda \delta_\lambda\left(Y_i\right)  \nonumber  
\end{align}
and hence
\begin{align}
\nabla_\gamma \circ \delta_\lambda=\lambda \delta_\lambda \nabla_\gamma.   \nonumber 
\end{align}
The respective distance function from the origin to some point in space $z=(x,y)\in\mathbb{R}^{n}$:
\begin{align}
\rho=\rho(z):=\left(|x|^{2(1+\gamma)}+(1+\gamma)^2|y|^2\right)^{\frac{1}{2(1+\gamma)}} . \nonumber
\end{align}
The function $\rho$ is positive, smooth away from the origin and symmetric. One can verify that
\begin{align*}
\nabla_{\gamma}\rho=\left(\frac{|x|^{2\gamma}x}{\rho^{2\gamma+1}},\frac{(1+\gamma)|x|^{\gamma}y}{\rho^{2\gamma+1}}\right),
\end{align*}
which gives us
\begin{align*}
|\nabla_{\gamma}\rho|=\frac{|x|^\gamma}{\rho^\gamma}.
\end{align*}
By direct calculation, we also get 
\begin{align}\label{main formula}
\nabla_{\gamma}\cdot\left(\rho^{c}|x|^{s}\nabla_{\gamma}\rho\right)=(Q+c+s-1)\frac{|x|^{2\gamma+s}}{\rho^{2\gamma+1-c}}
\end{align}
with $c,s\in\mathbb{R}$, which by setting $c=\gamma$ and $s=-\gamma$ gives
\begin{align}\label{formula for term2}
\nabla_{\gamma}\left(\frac{\rho^{\gamma}}{|x|^{\gamma}}\nabla_{\gamma}\rho\right)=\nabla_{\gamma}\left(\frac{\nabla_{\gamma}\rho}{|\nabla_{\gamma}\rho|}\right)=(Q-1)\frac{|x|^{\gamma}}{\rho^{\gamma+1}}.
\end{align}
We define $B_{R}=\{z\in\mathbb{R}^n:\rho(z)<R\}$ as a $\rho$-ball.

\section{Refined general weighted $L^p$-Hardy type inequalities and identities}\label{refined general}

In this section, we prove the improved general weighted $L^p$-Hardy type inequalities and identities related to the Baouendi-Grushin operator.

\begin{thm}\label{thm1}
Let $1<p<\infty$ and let $\Omega\subseteq\mathbb{R}^{m+k}$ be an open set such that the integrals below make sense. Let $v\in C^{1}(\Omega\backslash\Sigma)$ and $w\in L^{1}_{loc}(\Omega\backslash\Sigma)$ be nonnegative functions satisfying the condition (in the weak sense)
\begin{align}\label{main condt}
0\leq\phi:=\nabla_{\gamma}\cdot \left(w^{\frac{p-1}{p}}v^{\frac{1}{p}}\frac{\nabla_{\gamma}\rho}{|\nabla_{\gamma}\rho|}\right)-pw
\end{align}
a.e. in $\Omega\backslash\Sigma$ with $\Sigma$ being the set of singular points of $v$ and $w$.
\begin{enumerate}
\item Then, for all complex-valued $f\in C^{\infty}_0(\Omega\backslash\Sigma)$, there holds
\begin{align*}
\int_{\Omega}^{}v\frac{|\nabla_{\gamma}\rho\cdot\nabla_{\gamma}f|^{p}}{|\nabla_{\gamma}\rho|^{p}}dz\geq\int_{\Omega}^{}w|f|^{p}dz.
\end{align*}
    \item Furthermore, for all complex-valued $f\in C^{\infty}_0(\Omega\backslash\Sigma)$, we also have the identity
    \begin{align}\label{main res}
\int_{\Omega}^{}v\frac{|\nabla_{\gamma}\rho\cdot\nabla_{\gamma}f|^{p}}{|\nabla_{\gamma}\rho|^{p}}dz=\int_{\Omega}^{}w|f|^{p}dz+\int_{\Omega}^{}C_{p}(\xi,\eta)dz+\int_{\Omega}\phi|f|^{p}dz,
\end{align}
where the functional $C_p(\cdot,\cdot)$ and $\xi,\eta$ are given in (\ref{cp formula intro}) and (\ref{main notation intro}), respectively.
\end{enumerate} 
\end{thm}
\begin{flushleft}
\textit{Proof of Theorem \ref{thm1}}: Assume that (\ref{main condt}) holds, then
\end{flushleft}
\begin{align}\label{cancel term}
\int_{\Omega}^{}w|f|^{p}dz&= \frac{1}{p}\int_{\Omega}^{}|f|^{p}\text{div}_{\nabla_{\gamma}}\left(w^{\frac{p-1}{p}}v^{\frac{1}{p}}\frac{\nabla_{\gamma}\rho}{|\nabla_{\gamma}\rho|}\right)dz-\frac{1}{p}\int_{\Omega}\phi |f|^{p}dz \nonumber
\\&=-\text{Re}\int_{\Omega}^{}f|f|^{p-2}\frac{\overline{\nabla_{\gamma}\rho \cdot \nabla_{\gamma}f}}{|\nabla_{\gamma}\rho|}w^{\frac{p-1}{p}}v^{\frac{1}{p}}dz-\frac{1}{p}\int_{\Omega}\phi|f|^{p}dz.
\end{align}
By H{\"o}lder inequality, we get
\begin{align*}
\int_{\Omega}w|f|^{p}dz&\leq\int_{\Omega}|f|^{p-1}w^{\frac{p-1}{p}}v^{\frac{1}{p}}\frac{|\nabla_{\gamma}\rho\cdot\nabla_{\gamma}f|}{|\nabla_{\gamma}\rho|}dz \nonumber
\\&\leq\left(\int_{\Omega}w|f|^{p}dz\right)^{\frac{p-1}{p}}\left(\int_{\Omega}v\frac{|\nabla_{\gamma}\rho\cdot\nabla_{\gamma}f|^{p}}{|\nabla_{\gamma}\rho|^{p}}dz\right)^{\frac{1}{p}}.
\end{align*}
This immediately gives us \textit{Part (1)} of Theorem \ref{thm1}. Now for \textit{Part (2)}, we utilize the notation (\ref{main notation intro}) and formula (\ref{cp formula intro}):
\begin{align*}
&C_{p}(\xi,\eta)=v\frac{|\nabla_{\gamma}\rho\cdot\nabla_{\gamma}f|^{p}}{|\nabla_{\gamma}\rho|^{p}}-\left|-w^{\frac{1}{p}}f\right|^{p}-p\left|-w^{\frac{1}{p}}f\right|^{p-2}\text{Re}\left(-w^{\frac{1}{p}}f\right)
\\&\times\overline{\left(v^{\frac{1}{p}}\frac{\nabla_{\gamma}\rho \cdot \nabla_{\gamma}f}{|\nabla_{\gamma}\rho|} + w^{\frac{1}{p}}f\right)}.
\end{align*}
Expanding and simplifying results in
\begin{align*}
&C_{p}(\xi,\eta)=v\frac{|\nabla_{\gamma}\rho \cdot \nabla_{\gamma}f|^{p}}{|\nabla_{\gamma}\rho|^{p}}-w|f|^{p}+p\text{Re}w^{\frac{p-1}{p}}f|f|^{p-2}\overline{v^{\frac{1}{p}}\frac{\nabla_{\gamma}\rho \cdot \nabla_{\gamma}f}{|\nabla_{\gamma}\rho|}}+pw|f|^{p}.
\end{align*}
Integrating both sides, we have
\begin{align*}
&\int_{\Omega}^{}C_{p}(\xi,\eta)dz=\int_{\Omega}^{}v\frac{|\nabla_{\gamma}\rho\cdot\nabla_{\gamma}f|^{p}}{|\nabla_{\gamma}\rho|^{p}}dz-\int_{\Omega}^{}w|f|^{p}dz
\\&+p\text{Re}\int_{\Omega}^{}f|f|^{p-2}\frac{\overline{\nabla_{\gamma}\rho \cdot \nabla_{\gamma}f}}{|\nabla_{\gamma}\rho|}w^{\frac{p-1}{p}}v^{\frac{1}{p}}dz+p\int_{\Omega}^{}w|f|^{p}dz.
\end{align*}
Using the identity (\ref{cancel term}):
\begin{align*}
&\int_{\Omega}^{}C_{p}(\xi,\eta)dz=\int_{\Omega}^{}v\frac{|\nabla_{\gamma}\rho\cdot\nabla_{\gamma}f|^{p}}{|\nabla_{\gamma}\rho|^{p}}dz-\int_{\Omega}^{}w|f|^{p}dz
\\&+p\text{Re}\int_{\Omega}^{}f|f|^{p-2}\frac{\overline{\nabla_{\gamma}\rho \cdot \nabla_{\gamma}f}}{|\nabla_{\gamma}\rho|}w^{\frac{p-1}{p}}v^{\frac{1}{p}}dz+p\int_{\Omega}^{}w|f|^{p}dz
\\&=\int_{\Omega}^{}v\frac{|\nabla_{\gamma}\rho\cdot\nabla_{\gamma}f|^{p}}{|\nabla_{\gamma}\rho|^{p}}dz-\int_{\Omega}^{}w|f|^{p}dz
\\&-p\int_{\Omega}^{}w|f|^{p}dz-\int_{\Omega}\phi|f|^{p}dz+p\int_{\Omega}^{}w|f|^{p}dz
\\&=\int_{\Omega}^{}v\frac{|\nabla_{\gamma}\rho\cdot\nabla_{\gamma}f|^{p}}{|\nabla_{\gamma}\rho|^{p}}dz-\int_{\Omega}^{}w|f|^{p}dz-\int_{\Omega}\phi|f|^{p}dz.
\end{align*}
As a result, we obtain (\ref{main res}). \QEDA

\vspace{3mm}

If the condition (\ref{main condt}) is satisfied with $\phi=0$ in (\ref{main res}), then we obtain the remainder estimate from below for $p\geq2$ via the result from \cite{cazacu2024hardy}:

\begin{thm}\label{laptev thm}
Let $p\geq2$ and suppose the condition (\ref{main condt}) from Theorem \ref{thm1} is satisfied with $\phi=0$. That is, for all complex-valued $f\in C^{\infty}_0(\Omega\backslash\Sigma)$, we have
\begin{align}\label{laptev condt}
\int_{\Omega}^{}v\frac{|\nabla_{\gamma}\rho\cdot\nabla_{\gamma}f|^{p}}{|\nabla_{\gamma}\rho|^{p}}dz=\int_{\Omega}^{}w|f|^{p}dz+\int_{\Omega}^{}C_{p}(\xi,\eta)dz,
\end{align}
where the functional $C_p(\cdot,\cdot)$ and $\xi,\eta$ are given in (\ref{cp formula intro}) and (\ref{main notation intro}), respectively. Then, for all complex-valued $f\in C^{\infty}_0(\Omega\backslash\Sigma)$, we have
\begin{align}\label{laptev geq}
\int_{\Omega}C_{p}(\xi,\eta)\geq c_p\int_{\Omega}\left|v^{\frac{1}{p}}\frac{\nabla_{\gamma}\rho \cdot \nabla_{\gamma}f}{|\nabla_{\gamma}\rho|} + w^{\frac{1}{p}}f\right|^{p}dz,
\end{align}
where
\begin{align}\label{laptev cp const}
c_p
= \inf_{(s,t)\in\mathbb{R}^2\setminus\{(0,0)\}}
\frac{\bigl[t^2 + s^2 + 2s + 1\bigr]^{\frac p2} \;-\; 1 \;-\; p\,s}
{\bigl[t^2 + s^2\bigr]^{\frac p2}}
\;\in\;(0,1].
\end{align}
\end{thm}

Additionally, under the same condition $\phi=0$ in (\ref{main res}), recent results from \cite{CT24} allow us to obtain (\ref{main res}) with different remainder terms for $1<p<2\leq n$:

\begin{thm}\label{thm chi}
Let $1<p<2\leq n$ and suppose the condition (\ref{main condt}) from Theorem \ref{thm1} is satisfied with $\phi=0$. That is, for all complex-valued $f\in C^{\infty}_0(\Omega\backslash\Sigma)$, we have
\begin{align}\label{thm chi condt}
\int_{\Omega}^{}v\frac{|\nabla_{\gamma}\rho\cdot\nabla_{\gamma}f|^{p}}{|\nabla_{\gamma}\rho|^{p}}dz=\int_{\Omega}^{}w|f|^{p}dz+\int_{\Omega}^{}C_{p}(\xi,\eta)dz,
\end{align}
where the functional $C_p(\cdot,\cdot)$ and $\xi,\eta$ are given in (\ref{cp formula intro}) and (\ref{main notation intro}), respectively.
\begin{enumerate}
    \item  Then, for all complex-valued $f\in C^{\infty}_0(\Omega\backslash\Sigma)$, we have
    \begin{multline}\label{thm chi geq c1}
    \int_{\Omega}^{}C_{p}(\xi,\eta)dz\geq c_{1}(p)\int_{\Omega}\left(\left|v^{\frac{1}{p}}\frac{\nabla_{\gamma}\rho \cdot \nabla_{\gamma}f}{|\nabla_{\gamma}\rho|}\right|+|w^{\frac{1}{p}}f|\right)^{p-2}\biggl|v^{\frac{1}{p}}\frac{\nabla_{\gamma}\rho \cdot \nabla_{\gamma}f}{|\nabla_{\gamma}\rho|}\\ + w^{\frac{1}{p}}f\biggr|^{2}dz,
    \end{multline}
    where $c_1(p)$ is an explicit constant defined by
    \begin{align}\label{thm chi c1}
c_1(p) := \inf_{s^2 + t^2 > 0} \frac{\left( t^2 + s^2 + 2s + 1 \right)^{\frac{p}{2}} - 1 - ps}{\left( \sqrt{t^2 + s^2 + 2s + 1} + 1 \right)^{p-2} (t^2 + s^2)} \in \left( 0, \frac{p(p-1)}{2p-1} \right].
\end{align}
\item Moreover, for all complex-valued $f\in C^{\infty}_0(\Omega\backslash\Sigma)$, the remainder term is optimal since
\begin{align}\label{thm chi leq c2}
\int_{\Omega}^{}C_{p}(\xi,\eta)dz\leq c_{2}(p)\int_{\Omega}\left(\left|v^{\frac{1}{p}}\frac{\nabla_{\gamma}\rho \cdot \nabla_{\gamma}f}{|\nabla_{\gamma}\rho|}\right|+|w^{\frac{1}{p}}f|\right)^{p-2}\left|v^{\frac{1}{p}}\frac{\nabla_{\gamma}\rho \cdot \nabla_{\gamma}f}{|\nabla_{\gamma}\rho|} + w^{\frac{1}{p}}f\right|^{2}dz,
\end{align}
where $c_{2}(p)$ is an explicit constant defined by
\begin{align}\label{thm chi c2}
c_2(p) := \sup_{s^2 + t^2 > 0} \frac{\left( t^2 + s^2 + 2s + 1 \right)^{\frac{p}{2}} - 1 - ps}{\left( \sqrt{t^2 + s^2 + 2s + 1} + 1 \right)^{p-2} (t^2 + s^2)} \in \left[ \frac{p}{2^{p-1}}, +\infty \right).
\end{align}
\item In addition, for all complex-valued $f\in C^{\infty}_0(\Omega\backslash\Sigma)$, we have
\begin{multline}\label{thm chi geq c3}
\int_{\Omega}^{}C_{p}(\xi,\eta)dz\geq c_{3}(p)\int_{\Omega}\min\biggl\{\left|v^{\frac{1}{p}}\frac{\nabla_{\gamma}\rho \cdot \nabla_{\gamma}f}{|\nabla_{\gamma}\rho|} + w^{\frac{1}{p}}f\right|^{p},\\|w^{\frac{1}{p}}f|^{p-2}\left|v^{\frac{1}{p}}\frac{\nabla_{\gamma}\rho \cdot \nabla_{\gamma}f}{|\nabla_{\gamma}\rho|} + w^{\frac{1}{p}}f\right|^{2}\biggr\}dz,
\end{multline}
where $c_{3}(p)$ is an explicit constant defined by
\begin{multline}\label{thm chi c3}
c_3(p) := \min \biggl\{
\inf_{s^2 + t^2 \geq 1} \frac{(t^2 + s^2 + 2s + 1)^{\frac{p}{2}} - 1 - ps}{(t^2 + s^2)^{\frac{p}{2}}},\\
\inf_{0 < s^2 + t^2 < 1} \frac{(t^2 + s^2 + 2s + 1)^{\frac{p}{2}} - 1 - ps}{t^2 + s^2}
\biggr\}\in \left( 0, \frac{p(p-1)}{2} \right].
\end{multline}
\end{enumerate}
\end{thm}
Before proving Theorem \ref{laptev thm} and \ref{thm chi}, we first present the following lemmata from \cite{cazacu2024hardy} and \cite{CT24}, which play an important role in the argument:
\begin{lem}[\text{\cite[Step 3 of Proof of Lemma 3.4]{cazacu2024hardy}}]\label{lem4}
Let $p\geq2$. Then, for $\xi,\eta\in\mathbb{C}^n$, we have
\begin{align*}
C_{p}(\xi,\eta)\geq c_p|\eta|^{p},
\end{align*}
where $c_p$ is an explicit constant defined in (\ref{laptev cp const}).
\end{lem}
\begin{lem}[\text{\cite[Lemma 2.2]{CT24}}]
Let $1<p<2\leq n$. Then, for $\xi,\eta\in \mathbb{C}^{n}$, we have
\begin{align*}
C_p(\xi, \eta) \geq c_1(p) \left( |\xi| + |\xi - \eta| \right)^{p-2} |\eta|^2,
\end{align*}
where $c_{1}(p)$ is an explicit constant defined in (\ref{thm chi c1}).
\end{lem}
\begin{lem}[\text{\cite[Lemma 2.3]{CT24}}]
Let $1<p<2\leq n$. Then, for $\xi,\eta\in \mathbb{C}^{n}$, we have  
\begin{align*}
C_p(\xi, \eta) \leq c_2(p) \left( |\xi| + |\xi - \eta| \right)^{p-2} |\eta|^2,
\end{align*}
where $c_{2}(p)$ is an explicit constant defined in (\ref{thm chi c2}).
\end{lem}
\begin{lem}[\text{\cite[Lemma 2.4]{CT24}}]
Let $1<p<2\leq n$. Then, for $\xi,\eta\in \mathbb{C}^{n}$, we have  
\begin{align*}
C_p(\xi, \eta) \geq c_3(p) \min \left\{ |\eta|^p, |\xi - \eta|^{p-2} |\eta|^2 \right\}, 
\end{align*}
where $c_{3}(p)$ is an explicit constant defined in (\ref{thm chi c3}).
\end{lem}

\begin{flushleft}
\textit{Proof of Theorem \ref{laptev thm}}: Assume that (\ref{laptev condt}) holds, then we can write
\end{flushleft}
\begin{align*}
\int_{\Omega}^{}C_{p}(\xi,\eta)dz=\int_{\Omega}^{}v\frac{|\nabla_{\gamma}\rho\cdot\nabla_{\gamma}f|^{p}}{|\nabla_{\gamma}\rho|^{p}}dz-\int_{\Omega}^{}w|f|^{p}dz.
\end{align*}
Using Lemma \ref{lem4}, we derive (\ref{laptev geq}). \QEDA

\begin{flushleft}
\textit{Proof of Theorem \ref{thm chi}}: Assume that (\ref{thm chi condt}) holds, then, in the same way as in the proof of Theorem \ref{laptev thm}, we derive (\ref{thm chi geq c1}), (\ref{thm chi leq c2}) and (\ref{thm chi geq c3}). \QEDA
\end{flushleft}

\section{Applications of Theorem \ref{thm1}, \ref{laptev thm} and \ref{thm chi}}\label{applications}

Let $\epsilon>0$ be given. Define 
\begin{align*}
\rho_{\epsilon}:=\left(|x|^{2(1+\gamma)}_{\epsilon}+(1+\gamma)^{2}|y|^{2}\right)^{\frac{1}{2(1+\gamma)}},
\end{align*}
where $|x|_{\epsilon}:=\left(\epsilon^{2}+\sum_{i=1}^{m}x^{2}_{i}\right)^{\frac{1}{2}}$. To ensure the rigor of the following arguments, we need to replace $\rho$ with its regularization $\rho_{\epsilon}$ and take the limit $\epsilon\longrightarrow 0$ after performing the computations. However, for simplicity, we proceed with a formal approach. 

\subsection{Refinements of weighted $L^{p}$-Hardy type inequalities and identities.}\label{subsection1}

    \hspace{1mm}

\begin{flushleft}
We now apply Theorem \ref{thm1} to obtain refinements of previously known results on the whole space and $\rho$-ball $B_{R}$. For example, if we set    
\end{flushleft} 
\begin{align*}
v = 1, \quad w=\left(\frac{p-1}{p}\right)^p\frac{|x|^{\gamma p}}{(R-\rho)^{p}\rho^{\gamma p}}
\end{align*}
then, checking the condition (\ref{main condt}), we see that it actually holds and, thus, we have the following refined inequalities and identities of the result by Niu, Chen and Han \cite[Theorem 2.1]{NCH04}:
\begin{cor}\label{cor1}
Let $1<p<\infty$ and $R>0$. 
\begin{enumerate}
\item Then, for all complex-valued $f\in C^{\infty}_{0}(B_{R}\backslash\{(0,0)\})$, there holds
\begin{align}\label{cor1 eq0}
\int_{B_{R}}\frac{|\nabla_{\gamma}\rho\cdot\nabla_{\gamma}f|^p}{|\nabla_{\gamma}\rho|^p}dz\geq\left(\frac{p-1}{p}\right)^p\int_{B_R}\frac{|x|^{\gamma p}}{(R-\rho)^p\rho^{\gamma p}}|f|^pdz,
\end{align}
where the constant $\left(\frac{p-1}{p}\right)^{p}$ is sharp.
    \item Furthermore, for all complex-valued $f\in C^{\infty}_{0}(B_{R}\backslash\{(0,0)\})$, we also have the identity
\begin{multline}\label{cor1 eq}
\int_{B_{R}}\frac{|\nabla_{\gamma}\rho\cdot\nabla_{\gamma}f|^p}{|\nabla_{\gamma}\rho|^p}dz=\left(\frac{p-1}{p}\right)^p\int_{B_R}\frac{|x|^{\gamma p}}{(R-\rho)^p\rho^{\gamma p}}|f|^pdz\\+\int_{B_R}C_{p}\left(\frac{\nabla_{\gamma}\rho \cdot \nabla_{\gamma}f}{|\nabla_{\gamma}\rho|},\frac{\nabla_{\gamma}\rho \cdot \nabla_{\gamma}f}{|\nabla_{\gamma}\rho|} + \left(\frac{p-1}{p}\right)\frac{|x|^{\gamma }}{(R-\rho)\rho^{\gamma}}f\right)dz
\\+\left(\frac{p-1}{p}\right)^{p-1}(Q-1)\int_{B_{R}}\frac{|x|^{\gamma p}}{(R-\rho)^{p-1}\rho^{\gamma p+1}}|f|^{p}dz
\end{multline}
with functional $C_{p}(\cdot,\cdot)$ given in Theorem \ref{thm1}.
\end{enumerate}
\end{cor}
\begin{flushleft}
\textit{Proof of Corollary \ref{cor1}:} To obtain (\ref{cor1 eq0}) and (\ref{cor1 eq}), we need to check the condition (\ref{main condt}). If the condition is satisfied, then the desired results are obtained. First, we calculate the divergence:
\end{flushleft}
\begin{align}
&\nabla_{\gamma} \cdot\left(\left(\frac{p-1}{p}\right)^{p\frac{p-1}{p}}\left(\frac{|x|^{\gamma p}}{(R-\rho)^{p}\rho^{\gamma p}}\right)^{\frac{p-1}{p}}\frac{\nabla_{\gamma}\rho}{|\nabla_{\gamma}\rho|}\right)=\nabla_{\gamma}\cdot\left(\Phi \frac{\nabla_{\gamma}\rho}{|\nabla_{\gamma}\rho|}\right) \nonumber
\\&=\nabla_{\gamma}\Phi \cdot \frac{\nabla_{\gamma}\rho}{|\nabla_{\gamma}\rho|}+\Phi \nabla_{\gamma}\cdot \left(\frac{\nabla_{\gamma}\rho}{|\nabla_{\gamma}\rho|}\right) \nonumber
\\&=\frac{1}{|\nabla_{\gamma}\rho|}\nabla_{\gamma}\Phi \cdot \nabla_{\gamma}\rho + \Phi\nabla_{\gamma}\cdot \left(\frac{\nabla_{\gamma}\rho}{|\nabla_{\gamma}\rho|}\right) \nonumber
\\&=T_1+T_2,\label{cor1 t1+t2}
\end{align}
where
\begin{align}
&T_1=\frac{1}{|\nabla_{\gamma}\rho|}\nabla_{\gamma}\Phi \cdot \nabla_{\gamma}\rho, \nonumber
\\&T_2=\Phi\nabla_{\gamma}\cdot \left(\frac{\nabla_{\gamma}\rho}{|\nabla_{\gamma}\rho|}\right) \label{cor1 t2}
\end{align}
with
\begin{align*}
&\Phi=\left[\left(\frac{p-1}{p}\right)^{p}\frac{|x|^{\gamma p}}{(R-\rho)^{p}\rho^{\gamma p}}\right]^{\frac{p-1}{p}}
\\&=\left(\frac{p-1}{p}\right)^{p-1}\frac{|x|^{\gamma(p-1)}}{(R-\rho)^{p-1}\rho^{\gamma(p-1)}}.
\end{align*}
Substituting formula (\ref{formula for term2}) to (\ref{cor1 t2}), we get
\begin{align}
&T_2=\Phi \cdot\left[(Q-1)\frac{|x|^{\gamma}}{\rho^{\gamma+1}}\right] \nonumber
\\&=\left(\frac{p-1}{p}\right)^{p-1}\frac{|x|^{\gamma(p-1)}}{(R-\rho)^{p-1}\rho^{\gamma(p-1)}}\cdot(Q-1)\frac{|x|^{\gamma}}{\rho^{\gamma+1}} \nonumber
\\&=\left(\frac{p-1}{p}\right)^{p-1}(Q-1)\frac{|x|^{\gamma p}}{(R-\rho)^{p-1}\rho^{\gamma p+1}}. \label{cor1 t2_1}
\end{align}
Now we compute $T_1$:
\begin{align}
&T_1=\frac{1}{|\nabla_{\gamma}\rho|}(\nabla_{\gamma}\Phi)\cdot \nabla_{\gamma}\rho \nonumber
\\&=\frac{1}{|\nabla_{\gamma}\rho|}\left(\frac{p-1}{p}\right)^{p-1}\nabla_{\gamma}\left((R-\rho)^{-(p-1)}\left(\frac{|x|^{\gamma}}{\rho^{\gamma}}\right)^{p-1}\right)\cdot\nabla_{\gamma}\rho \nonumber
\\&=\frac{\rho^{\gamma}}{|x|^{\gamma}}\left(\frac{p-1}{p}\right)^{p-1}(p-1)(R-\rho)^{-p}\frac{|x|^{\gamma(p+1)}}{\rho^{\gamma(p+1)}} \nonumber
\\&=p\left(\frac{p-1}{p}\right)^{p}\frac{|x|^{\gamma p}}{(R-\rho)^{p}\rho^{\gamma p}}. \label{cor1 t1_1}
\end{align}
Putting (\ref{cor1 t2_1}) and (\ref{cor1 t1_1}) to (\ref{cor1 t1+t2}), we obtain
\begin{align*}
&T_{1}+T_{2}=p\left(\frac{p-1}{p}\right)^{p}\frac{|x|^{\gamma p}}{(R-\rho)^{p}\rho^{\gamma p}}+\left(\frac{p-1}{p}\right)^{p-1}(Q-1)\frac{|x|^{\gamma p}}{(R-\rho)^{p-1}\rho^{\gamma p+1}}
\\&=pw+\left(\frac{p-1}{p}\right)^{p-1}(Q-1)\frac{|x|^{\gamma p}}{(R-\rho)^{p-1}\rho^{\gamma p+1}}.
\end{align*}
This implies that $\phi=\left(\frac{p-1}{p}\right)^{p-1}(Q-1)\frac{|x|^{\gamma p}}{(R-\rho)^{p-1}\rho^{\gamma p+1}}\geq0$. Therefore, the condition (\ref{main condt}) is satisfied and the results (\ref{cor1 eq0}) and (\ref{cor1 eq}) are obtained. Since the function $h=\left(R-\rho\right)^{-\frac{p-1}{p}}$ satisfies the H{\"o}lder equality condition
\begin{align*}
\left(\frac{p}{p-1}\right)^{p}\frac{|\nabla_{\gamma}\rho\cdot\nabla_{\gamma}h|^{p}}{|\nabla_{\gamma}\rho|^{p}}=\frac{|x|^{\gamma p}}{(R-\rho)^{p}\rho^{\gamma p}}|h|^{p},
\end{align*}
the constant in the inequality (\ref{cor1 eq0}) is sharp. \QEDA

\vspace{2mm}

We are able to obtain the refinements of some results by D'Ambrosio in \cite{D'A04} and \cite{D'A05}. For instance, setting 
\begin{align*}
v=|x|^{\beta-\gamma p}\rho^{p(1+\gamma)-\alpha},\quad w=\left(\frac{Q+\beta-\alpha}{p}\right)^{p}\frac{|x|^{\beta}}{\rho^{\alpha}},
\end{align*}
we get the following refinements of the weighted $L^p$-Hardy type inequality \cite[Formula (3.1) of Theorem 3.1]{D'A04}:
\begin{cor}\label{cor2}
Let $1<p<\infty$, $m,k\geq1$ and $\alpha,\beta\in \mathbb{R}$ be such that $Q> \alpha-\beta$.  
\begin{enumerate}
\item Then, for all complex-valued $f\in C^{\infty}_{0}(\mathbb{R}^m\times\mathbb{R}^{n}\backslash\{(0,0)\})$, there holds
\begin{align}\label{cor2 eq0}
\int_{\mathbb{R}^n}|x|^{\beta-\gamma p}\rho^{p(1+\gamma)-\alpha}\frac{|\nabla_{\gamma}\rho\cdot\nabla_{\gamma}f|^p}{|\nabla_{\gamma}\rho|^p}dz\geq\left(\frac{Q+\beta-\alpha}{p}\right)^{p}\int_{\mathbb{R}^n}\frac{|x|^{\beta}}{\rho^{\alpha}}|f|^pdz,
\end{align}
where the constant $\left(\frac{Q+\beta-\alpha}{p}\right)^{p}$ is sharp.
    \item Furthermore, for all complex-valued $f\in C^{\infty}_{0}(\mathbb{R}^m\times\mathbb{R}^{n}\backslash\{(0,0)\})$, we also have the identity
\begin{multline}\label{cor2 eq}
\int_{\mathbb{R}^n}|x|^{\beta-\gamma p}\rho^{p(1+\gamma)-\alpha}\frac{|\nabla_{\gamma}\rho\cdot\nabla_{\gamma}f|^p}{|\nabla_{\gamma}\rho|^p}dz=\left(\frac{Q+\beta-\alpha}{p}\right)^{p}\int_{\mathbb{R}^n}\frac{|x|^{\beta}}{\rho^{\alpha}}|f|^pdz
\\+\int_{\mathbb{R}^n}C_{p}\biggl(|x|^{\frac{\beta-\gamma p}{p}}\rho^{\frac{p(1+\gamma)-\alpha}{p}}\frac{\nabla_{\gamma}\rho \cdot \nabla_{\gamma}f}{|\nabla_{\gamma}\rho|},|x|^{\frac{\beta-\gamma p}{p}}\rho^{\frac{p(1+\gamma)-\alpha}{p}}\frac{\nabla_{\gamma}\rho \cdot \nabla_{\gamma}f}{|\nabla_{\gamma}\rho|} \\+ \left(\frac{Q+\beta-\alpha}{p}\right)\frac{|x|^{\frac{\beta}{p}}}{\rho^{\frac{\alpha}{p}}}f\biggr)dz
\end{multline}
with functional $C_{p}(\cdot,\cdot)$ given in Theorem \ref{thm1}.
\end{enumerate}
\end{cor}

Substituting $\beta=\gamma p$ and $\alpha=p(1+\gamma)$ in (\ref{cor2 eq}), we obtain the following identity:

\begin{cor}
Let $1<p<\infty$ and $\gamma\geq0$. Then, we get the following identity:
\begin{multline}\label{gar ref}
\int_{\mathbb{R}^n}\frac{|\nabla_{\gamma}\rho\cdot\nabla_{\gamma}f|^{p}}{|\nabla_{\gamma}\rho|^{p}}dz=
\left(\frac{Q-p}{p}\right)^{p}\int_{\mathbb{R}^n}\frac{|x|^{\gamma p}}{\rho^{\gamma p+p}}|f|^{p}dz
\\+\int_{\mathbb{R}^n}C_{p}\left(\frac{\nabla_{\gamma}\rho \cdot \nabla_{\gamma}f}{|\nabla_{\gamma}\rho|},\frac{\nabla_{\gamma}\rho \cdot \nabla_{\gamma}f}{|\nabla_{\gamma}\rho|}+\left(\frac{Q-p}{p}\right)\frac{|x|^{\gamma}}{\rho^{\gamma+1}}f\right)dz
\end{multline}
for any complex-valued $f\in C^{\infty}_{0}(\mathbb{R}^m\times\mathbb{R}^{n}\backslash\{(0,0)\})$.
\end{cor}
\begin{rem}
By setting \( \gamma = 0 \) in (\ref{gar ref}), we obtain a refinement of the classical Hardy inequality (\ref{class_multdim_hardy}). This result was also derived, in the special case, by Kalaman and the first author of this paper in \cite{kalaman2024cylindrical} under specific parameter choices: let $1<p<\infty$. Then, for any complex-valued function \( f \in C^{\infty}_{0}(\mathbb{R}^m \times \mathbb{R}^n \setminus \{(0,0)\}) \), the following identity holds:
\begin{multline*}
\int_{\mathbb{R}^n} \frac{|z \cdot \nabla f|^p}{|z|^p} \, dz 
= \left( \frac{n - p}{p} \right)^p \int_{\mathbb{R}^n} \frac{|f|^p}{|z|^p} \, dz 
\\+ \int_{\mathbb{R}^n} C_p \left( \frac{z \cdot \nabla f}{|z|}, \frac{z \cdot \nabla f}{|z|} + \left( \frac{n - p}{p} \right) f \right) \, dz.
\end{multline*}
\end{rem}

\begin{rem}
In the case when $p=2$ in (\ref{gar ref}), we get the following refinement of (\ref{Gar}): let $\gamma\geq0$. Then, for all complex-valued \( f \in C^{\infty}_{0}(\mathbb{R}^m \times \mathbb{R}^n \setminus \{(0,0)\}) \), we have
\begin{multline*}
\int_{\mathbb{R}^n}\frac{|\nabla_{\gamma}\rho\cdot\nabla_{\gamma}f|^{2}}{|\nabla_{\gamma}\rho|^{2}}dz=
\left(\frac{Q-2}{2}\right)^{2}\int_{\mathbb{R}^n}\frac{|x|^{2\gamma}}{\rho^{2\gamma+2}}|f|^{2}dz
\\+\int_{\mathbb{R}^n}C_{2}\left(\frac{\nabla_{\gamma}\rho \cdot \nabla_{\gamma}f}{|\nabla_{\gamma}\rho|},\frac{\nabla_{\gamma}\rho \cdot \nabla_{\gamma}f}{|\nabla_{\gamma}\rho|}+\left(\frac{Q-2}{2}\right)\frac{|x|^{\gamma}}{\rho^{\gamma+1}}f\right)dz.
\end{multline*}
\end{rem}

Since Corollary \ref{cor2} gives $\phi=0$ in (\ref{cor2 eq}), we can apply Theorem \ref{laptev thm} to obtain remainder estimate of (\ref{cor2 eq}) from below for $p\geq2$:

\begin{cor}\label{dam laptev}
Let $p\geq2$, $m,k\geq1$ and $\alpha,\beta\in\mathbb{R}$ be such that $Q>\alpha-\beta$. Then, for any complex-valued $f\in C^{\infty}_{0}(\mathbb{R}^m\times\mathbb{R}^{n}\backslash\{(0,0)\})$, we have
\begin{multline*}
\int_{\mathbb{R}^n}|x|^{\beta-\gamma p}\rho^{p(1+\gamma)-\alpha}\frac{|\nabla_{\gamma}\rho\cdot\nabla_{\gamma}f|^p}{|\nabla_{\gamma}\rho|^p}dz
    -\left(\frac{Q+\beta-\alpha}{p}\right)^{p}\int_{\mathbb{R}^n}\frac{|x|^{\beta}}{\rho^{\alpha}}|f|^pdz \\\geq c_p\int_{\mathbb{R}^n}\left||x|^{\frac{\beta-\gamma p}{p}}\rho^{\frac{p(1+\gamma)-\alpha}{p}}\frac{\nabla_{\gamma}\rho \cdot \nabla_{\gamma}f}{|\nabla_{\gamma}\rho|}+ \left(\frac{Q+\beta-\alpha}{p}\right)\frac{|x|^{\frac{\beta}{p}}}{\rho^{\frac{\alpha}{p}}}f\right|^{p}dz,
\end{multline*}
where $c_p$ is an explicit constant defined in (\ref{laptev cp const}).
\end{cor}

Furthermore, applying Theorem \ref{thm chi}, we obtain (\ref{cor2 eq}) with different remainder terms for $1<p<2\leq n$:

\begin{cor}\label{dam chinese}
Let $1<p<2\leq n$, $m,k\geq1$ and $\alpha,\beta\in \mathbb{R}$ be such that $Q>\alpha-\beta$. Then, for any complex-valued $f\in C^{\infty}_{0}(\mathbb{R}^m\times\mathbb{R}^{n}\backslash\{(0,0)\})$,
\begin{enumerate}
    \item for constants $c_{1}(p), c_{2}(p)>0$ defined in (\ref{thm chi c1}) and (\ref{thm chi c2}), respectively, we have
    \begin{align*}
    &c_{2}(p)\int_{\mathbb{R}^n}\Biggl(\left||x|^{\frac{\beta-\gamma p}{p}}\rho^{\frac{p(1+\gamma)-\alpha}{p}}\frac{\nabla_{\gamma}\rho \cdot \nabla_{\gamma}f}{|\nabla_{\gamma}\rho|}\right|+\left|\left(\frac{Q+\beta-\alpha}{p}\right)\frac{|x|^{\frac{\beta}{p}}}{\rho^{\frac{\alpha}{p}}}f\right|\Biggr)^{p-2} \nonumber
    \\&\times\left||x|^{\frac{\beta-\gamma p}{p}}\rho^{\frac{p(1+\gamma)-\alpha}{p}}\frac{\nabla_{\gamma}\rho \cdot \nabla_{\gamma}f}{|\nabla_{\gamma}\rho|}+ \left(\frac{Q+\beta-\alpha}{p}\right)\frac{|x|^{\frac{\beta}{p}}}{\rho^{\frac{\alpha}{p}}}f\right|^{2}dz \nonumber
    \\&\geq \int_{\mathbb{R}^n}|x|^{\beta-\gamma p}\rho^{p(1+\gamma)-\alpha}\frac{|\nabla_{\gamma}\rho\cdot\nabla_{\gamma}f|^p}{|\nabla_{\gamma}\rho|^p}dz
    -\left(\frac{Q+\beta-\alpha}{p}\right)^{p}\int_{\mathbb{R}^n}\frac{|x|^{\beta}}{\rho^{\alpha}}|f|^pdz \nonumber
    \\&\geq c_{1}(p)\int_{\mathbb{R}^n}\Biggl(\left||x|^{\frac{\beta-\gamma p}{p}}\rho^{\frac{p(1+\gamma)-\alpha}{p}}\frac{\nabla_{\gamma}\rho \cdot \nabla_{\gamma}f}{|\nabla_{\gamma}\rho|}\right|+\left|\left(\frac{Q+\beta-\alpha}{p}\right)\frac{|x|^{\frac{\beta}{p}}}{\rho^{\frac{\alpha}{p}}}f\right|\Biggr)^{p-2} \nonumber
    \\&\times\left||x|^{\frac{\beta-\gamma p}{p}}\rho^{\frac{p(1+\gamma)-\alpha}{p}}\frac{\nabla_{\gamma}\rho \cdot \nabla_{\gamma}f}{|\nabla_{\gamma}\rho|} + \left(\frac{Q+\beta-\alpha}{p}\right)\frac{|x|^{\frac{\beta}{p}}}{\rho^{\frac{\alpha}{p}}}f\right|^{2}dz;
    \end{align*}
    \item for constant $c_{3}(p)$ defined in (\ref{thm chi c3}), we have
    \begin{multline*}
    \int_{\mathbb{R}^n}|x|^{\beta-\gamma p}\rho^{p(1+\gamma)-\alpha}\frac{|\nabla_{\gamma}\rho\cdot\nabla_{\gamma}f|^p}{|\nabla_{\gamma}\rho|^p}dz
    -\left(\frac{Q+\beta-\alpha}{p}\right)^{p}\int_{\mathbb{R}^n}\frac{|x|^{\beta}}{\rho^{\alpha}}|f|^pdz 
    \\\geq c_{3}(p)\int_{\mathbb{R}^n}\min\Biggl\{\left||x|^{\frac{\beta-\gamma p}{p}}\rho^{\frac{p(1+\gamma)-\alpha}{p}}\frac{\nabla_{\gamma}\rho \cdot \nabla_{\gamma}f}{|\nabla_{\gamma}\rho|} + \left(\frac{Q+\beta-\alpha}{p}\right)\frac{|x|^{\frac{\beta}{p}}}{\rho^{\frac{\alpha}{p}}}f\right|^{p},\\\left|\left(\frac{Q+\beta-\alpha}{p}\right)\frac{|x|^{\frac{\beta}{p}}}{\rho^{\frac{\alpha}{p}}}f\right|^{p-2}\Biggl||x|^{\frac{\beta-\gamma p}{p}}\rho^{\frac{p(1+\gamma)-\alpha}{p}}\frac{\nabla_{\gamma}\rho \cdot \nabla_{\gamma}f}{|\nabla_{\gamma}\rho|}\\ + \left(\frac{Q+\beta-\alpha}{p}\right)\frac{|x|^{\frac{\beta}{p}}}{\rho^{\frac{\alpha}{p}}}f\Biggr|^{2}\Biggr\}dz.
    \end{multline*}
\end{enumerate}
\end{cor}

\begin{flushleft}
\textit{Proof of Corollary \ref{cor2}:} Calculating divergence:
\end{flushleft}
\begin{align}\label{damfor}
&\nabla_{\gamma} \cdot\left(\left(\left(\frac{Q+\beta-\alpha}{p}\right)^{p}\frac{|x|^{\beta }}{\rho^{\alpha}}\right)^{\frac{p-1}{p}}|x|^{\frac{\beta-\gamma p}{p}}\rho^{\frac{p(1+\gamma)-\alpha}{p}}\frac{\nabla_{\gamma}\rho}{|\nabla_{\gamma}\rho|}\right)\nonumber
\\&=\nabla_{\gamma}\cdot\left(\left(\frac{Q+\beta-\alpha}{p}\right)^{p-1}|x|^{\beta-\gamma}\rho^{1+\gamma-\alpha}\frac{\nabla_{\gamma}\rho}{|\nabla_{\gamma}\rho|}\right)\nonumber
\\&=\left(\frac{Q+\beta-\alpha}{p}\right)^{p-1}\nabla_{\gamma}\cdot\left(|x|^{\beta-2\gamma}\rho^{1+2\gamma-\alpha}\nabla_{\gamma}\rho\right).
\end{align}
Using formula (\ref{main formula}), in (\ref{damfor}), with $s=\beta-2\gamma$ and $c=1+2\gamma-\alpha$, we get
\begin{align*}
&\left(\frac{Q+\beta-\alpha}{p}\right)^{p-1}\nabla_{\gamma}\cdot\left(|x|^{\beta-2\gamma}\rho^{1+2\gamma-\alpha}\nabla_{\gamma}\rho\right)=\left(\frac{Q+\beta-\alpha}{p}\right)^{p-1}(Q+\beta-\alpha)\frac{|x|^{\beta}}{\rho^{\alpha}}
\\&=p\left(\frac{Q+\beta-\alpha}{p}\right)^{p}\frac{|x|^{\beta}}{\rho^{\alpha}}=pw+0.
\end{align*}
As we can see, the condition (\ref{main condt}) is satisfied with $\phi=0$, giving the desired results (\ref{cor2 eq0}) and (\ref{cor2 eq}). Since the function $h=\rho^{\left(\frac{Q+\beta-\alpha}{p}\right)}$ satisfies the H{\"o}lder equality condition,
\begin{align*}
\left(\frac{p}{Q+\beta-\alpha}\right)^{p}|x|^{\beta-\gamma p}\rho^{(1+\gamma)p-\alpha}\frac{|\nabla_{\gamma}h\cdot\nabla_{\gamma}\rho|^{p}}{|\nabla_{\gamma}\rho|^{p}}=\frac{|x|^{\beta}}{\rho^{\alpha}}|h|^{p},
\end{align*}
the constant in the inequality (\ref{cor2 eq0}) is sharp. \QEDA

\vspace{2mm}

Setting 

\begin{align*}
v = \frac{|\nabla_\gamma \rho|^{\alpha}}{\rho^{p(\theta - 1)}}, \quad w = \left(\frac{Q - p\theta}{p}\right)^{p}\frac{\left| \nabla_\gamma \rho \right|^{\alpha + p}}{\rho^{p\theta}}
\end{align*}
gives us an extended version of D'Arca results \cite[Theorem 4.5 and 4.6]{d2024weighted} for $1<p<\infty$, $\alpha\in\mathbb{R}$ and for complex-valued $f\in C^{\infty}_{0}(B_{R}\backslash\{(0,0)\})$:
\begin{cor}\label{cor4.5}
Let $1<p<\infty$, $\alpha,\theta\in\mathbb{R}$ and $R>0$ be such that $Q> p\theta$.
\begin{enumerate}
    \item Then, for all complex-valued $f\in C^{\infty}_{0}(B_{R}\backslash\{(0,0)\})$, there holds
    \begin{align}\label{cor4.5 eq0}
    \int_{B_{R}}\frac{|\nabla_\gamma \rho|^{\alpha}}{\rho^{p(\theta - 1)}}\frac{|\nabla_{\gamma}\rho\cdot\nabla_{\gamma}f|^p}{|\nabla_{\gamma}\rho|^p}dz\geq\left(\frac{Q-p\theta}{p}\right)^{p}\int_{B_{R}}\frac{\left| \nabla_\gamma \rho \right|^{\alpha + p}}{\rho^{p\theta}}|f|^{p}dz,
    \end{align}
    where the constant $\left(\frac{Q-p\theta}{p}\right)^{p}$ is sharp.
    \item Furthermore, for all complex-valued $f$ $\in$  $C^{\infty}_{0}(B_{R}\backslash\{(0,0)\})$, we also have the identity
\begin{multline}\label{cor4.5 eq}
\int_{B_{R}}\frac{|\nabla_\gamma \rho|^{\alpha}}{\rho^{p(\theta - 1)}}\frac{|\nabla_{\gamma}\rho\cdot\nabla_{\gamma}f|^p}{|\nabla_{\gamma}\rho|^p}dz=\left(\frac{Q-p\theta}{p}\right)^{p}\int_{B_{R}}\frac{\left| \nabla_\gamma \rho \right|^{\alpha + p}}{\rho^{p\theta}}|f|^{p}dz
\\+\int_{B_{R}}C_{p}\Biggl(\frac{|\nabla_\gamma \rho|^{\frac{\alpha}{p}}}{\rho^{(\theta - 1)}}\frac{\nabla_{\gamma}\rho \cdot \nabla_{\gamma}f}{|\nabla_{\gamma}\rho|}, \frac{|\nabla_\gamma \rho|^{\frac{\alpha}{p}}}{\rho^{(\theta - 1)}}\frac{\nabla_{\gamma}\rho \cdot \nabla_{\gamma}f}{|\nabla_{\gamma}\rho|} \\+ \left(\frac{Q - p\theta}{p}\right)\frac{\left| \nabla_\gamma \rho \right|^{\frac{\alpha + p}{p}}}{\rho^{\theta}}f\Biggr)dz
\end{multline}
with functional $C_{p}(\cdot,\cdot)$ given in Theorem \ref{thm1}.
\end{enumerate}
\end{cor}
Under $\theta=1$ and $\alpha=0$ in (\ref{cor4.5 eq}), we recover the result (\ref{gar ref}) on the $\rho$-ball $B_{R}$:
\begin{cor}
Let $1<p<\infty$ and $\gamma\geq0$. Then, we get the following identity:
\begin{multline*}
\int_{B_{R}}\frac{|\nabla_{\gamma}\rho\cdot\nabla_{\gamma}f|^{p}}{|\nabla_{\gamma}\rho|^{p}}dz=
\left(\frac{Q-p}{p}\right)^{p}\int_{B_{R}}\frac{|x|^{\gamma p}}{\rho^{\gamma p+p}}|f|^{p}dz
\\+\int_{B_{R}}C_{p}\left(\frac{\nabla_{\gamma}\rho \cdot \nabla_{\gamma}f}{|\nabla_{\gamma}\rho|},\frac{\nabla_{\gamma}\rho \cdot \nabla_{\gamma}f}{|\nabla_{\gamma}\rho|}+\left(\frac{Q-p}{p}\right)\frac{|x|^{\gamma}}{\rho^{\gamma+1}}f\right)dz
\end{multline*}
for any complex-valued $f\in C^{\infty}_{0}(B_{R}\backslash\{(0,0)\})$.
\end{cor}

In the same way as with (\ref{cor2 eq}), we have $\phi=0$ in (\ref{cor4.5 eq}), which, by applying Theorem \ref{laptev thm}, gives the following remainder estimate for $p\geq2$:

\begin{cor}
Let $p\geq2$ and $\alpha, \theta\in\mathbb{R}$ be such that $Q>p\theta$. Then, for all complex-valued $f$ $\in$  $C^{\infty}_{0}(B_{R}\backslash\{(0,0)\})$, we have
\begin{multline*}
\int_{B_{R}}\frac{|\nabla_\gamma \rho|^{\alpha}}{\rho^{p(\theta - 1)}}\frac{|\nabla_{\gamma}\rho\cdot\nabla_{\gamma}f|^p}{|\nabla_{\gamma}\rho|^p}dz-\left(\frac{Q-p\theta}{p}\right)^{p}\int_{B_{R}}\frac{\left| \nabla_\gamma \rho \right|^{\alpha + p}}{\rho^{p\theta}}|f|^{p}dz\\\geq c_p\int_{B_{R}}\left|\frac{|\nabla_\gamma \rho|^{\frac{\alpha}{p}}}{\rho^{\theta - 1}}\frac{\nabla_{\gamma}\rho \cdot \nabla_{\gamma}f}{|\nabla_{\gamma}\rho|} + \left(\frac{Q - p\theta}{p}\right)\frac{\left| \nabla_\gamma \rho \right|^{\frac{\alpha + p}{p}}}{\rho^{\theta}}f\right|^{p}dz,
\end{multline*}
where $c_p$ is an explicit constant defined in (\ref{laptev cp const}).
\end{cor}

Here, using Theorem \ref{thm chi}, we also obtain (\ref{cor4.5 eq}) with different remainder terms for $1<p<2\leq n$:

\begin{cor}
Let $1<p<2\leq n$ and $\alpha,\theta\in\mathbb{R}$ be such that $Q>p\theta$. Then, for any complex-valued $f\in C^{\infty}_{0}(B_{R}\backslash\{(0,0)\})$,
\begin{enumerate}
    \item for constants $c_{1}(p), c_{2}(p)>0$ defined in (\ref{thm chi c1}) and (\ref{thm chi c2}), respectively, we have
    \begin{align*}
    &c_{2}(p)\int_{B_{R}}\left(\left|\frac{|\nabla_\gamma \rho|^{\frac{\alpha}{p}}}{\rho^{\theta - 1}}\frac{\nabla_{\gamma}\rho \cdot \nabla_{\gamma}f}{|\nabla_{\gamma}\rho|}\right| + \left|\left(\frac{Q - p\theta}{p}\right)\frac{\left| \nabla_\gamma \rho \right|^{\frac{\alpha + p}{p}}}{\rho^{\theta}}f\right|\right)^{p-2}
    \\&\times\left|\frac{|\nabla_\gamma \rho|^{\frac{\alpha}{p}}}{\rho^{\theta - 1}}\frac{\nabla_{\gamma}\rho \cdot \nabla_{\gamma}f}{|\nabla_{\gamma}\rho|} + \left(\frac{Q - p\theta}{p}\right)\frac{\left| \nabla_\gamma \rho \right|^{\frac{\alpha + p}{p}}}{\rho^{\theta}}f\right|^{2}dz
    \\&\geq\int_{B_{R}}\frac{|\nabla_\gamma \rho|^{\alpha}}{\rho^{p(\theta - 1)}}\frac{|\nabla_{\gamma}\rho\cdot\nabla_{\gamma}f|^p}{|\nabla_{\gamma}\rho|^p}dz-\left(\frac{Q-p\theta}{p}\right)^{p}\int_{B_{R}}\frac{\left| \nabla_\gamma \rho \right|^{\alpha + p}}{\rho^{p\theta}}|f|^{p}dz
    \\&\geq c_{1}(p)\int_{B_{R}}\left(\left|\frac{|\nabla_\gamma \rho|^{\frac{\alpha}{p}}}{\rho^{\theta - 1}}\frac{\nabla_{\gamma}\rho \cdot \nabla_{\gamma}f}{|\nabla_{\gamma}\rho|}\right| + \left|\left(\frac{Q - p\theta}{p}\right)\frac{\left| \nabla_\gamma \rho \right|^{\frac{\alpha + p}{p}}}{\rho^{\theta}}f\right|\right)^{p-2}
    \\&\times\left|\frac{|\nabla_\gamma \rho|^{\frac{\alpha}{p}}}{\rho^{\theta - 1}}\frac{\nabla_{\gamma}\rho \cdot \nabla_{\gamma}f}{|\nabla_{\gamma}\rho|} + \left(\frac{Q - p\theta}{p}\right)\frac{\left| \nabla_\gamma \rho \right|^{\frac{\alpha + p}{p}}}{\rho^{\theta}}f\right|^{2}dz;
    \end{align*}
    \item for constant $c_{3}(p)$ defined in (\ref{thm chi c3}), we have
    \begin{multline*}
    \int_{B_{R}}\frac{|\nabla_\gamma \rho|^{\alpha}}{\rho^{p(\theta - 1)}}\frac{|\nabla_{\gamma}\rho\cdot\nabla_{\gamma}f|^p}{|\nabla_{\gamma}\rho|^p}dz-\left(\frac{Q-p\theta}{p}\right)^{p}\int_{B_{R}}\frac{\left| \nabla_\gamma \rho \right|^{\alpha + p}}{\rho^{p\theta}}|f|^{p}dz\\\geq c_{3}(p)\int_{B_{R}}\min\Biggl\{\left|\frac{|\nabla_\gamma \rho|^{\frac{\alpha}{p}}}{\rho^{\theta - 1}}\frac{\nabla_{\gamma}\rho \cdot \nabla_{\gamma}f}{|\nabla_{\gamma}\rho|} + \left(\frac{Q - p\theta}{p}\right)\frac{\left| \nabla_\gamma \rho \right|^{\frac{\alpha + p}{p}}}{\rho^{\theta}}f\right|^{p},\\\left|\left(\frac{Q - p\theta}{p}\right)\frac{\left| \nabla_\gamma \rho \right|^{\frac{\alpha + p}{p}}}{\rho^{\theta}}f\right|^{p-2}\left|\frac{|\nabla_\gamma \rho|^{\frac{\alpha}{p}}}{\rho^{\theta - 1}}\frac{\nabla_{\gamma}\rho \cdot \nabla_{\gamma}f}{|\nabla_{\gamma}\rho|} + \left(\frac{Q - p\theta}{p}\right)\frac{\left| \nabla_\gamma \rho \right|^{\frac{\alpha + p}{p}}}{\rho^{\theta}}f\right|^{2}\Biggr\}dz.
    \end{multline*}
\end{enumerate}
\end{cor}
\begin{flushleft}
\textit{Proof of Corollary \ref{cor4.5}:} Calculating divergence:
\end{flushleft}
\begin{align}
&\nabla_{\gamma}\cdot\left(\left(\left(\frac{Q-p\theta}{p}\right)^{p}\frac{|\nabla_{\gamma}\rho|^{\alpha+p}}{\rho^{p\theta}}\right)^{\frac{p-1}{p}}\left(\frac{|\nabla_{\gamma}\rho|^{\alpha}}{\rho^{p(\theta-1)}}\right)^{\frac{1}{p}}\frac{\nabla_{\gamma}\rho}{|\nabla_{\gamma}\rho|}\right) \nonumber
\\& = \nabla_{\gamma}\cdot\left(\left(\frac{Q-p\theta}{p}\right)^{p-1}\frac{|\nabla_{\gamma}\rho|^{\frac{(\alpha+p)(p-1)}{p}}}{\rho^{\theta(p-1)}}\frac{|\nabla_{\gamma}\rho|^{\frac{\alpha}{p}}}{\rho^{(\theta-1)}}\frac{\nabla_{\gamma}\rho}{|\nabla_{\gamma}\rho|}\right)   \nonumber
\\& = \left(\frac{Q-p\theta}{p}\right)^{p-1}\nabla_{\gamma}\cdot\left(\frac{|\nabla_{\gamma}\rho|^{\alpha+p-2}}{\rho^{p\theta-1}}\nabla_{\gamma}\rho\right) \nonumber
\\& = \left(\frac{Q-p\theta}{p}\right)^{p-1}\nabla_{\gamma}\cdot\left(\frac{|x|^{{\gamma}(\alpha+p-2)}}{\rho^{\gamma(\alpha+p-2)+p\theta-1}}\nabla_{\gamma}\rho\right)   \nonumber
\\& = \left(\frac{Q-p\theta}{p}\right)^{p-1}\nabla_{\gamma}\cdot\left(|x|^{\gamma(\alpha+p-2)}\rho^{-\gamma(\alpha+p-2)-p\theta+1}\nabla_{\gamma}\rho\right) .\label{darca proof}
\end{align}
Applying formula (\ref{main formula}) in (\ref{darca proof}) with $s=\gamma(\alpha+p-2)$ and $c=-\gamma(\alpha+p-2)-p\theta+1$, we obtain
\begin{align*}
&\left(\frac{Q-p\theta}{p}\right)^{p-1}\nabla_{\gamma}\cdot\left(|x|^{\gamma(\alpha+p-2)}\rho^{-\gamma(\alpha+p-2)-p\theta+1}\nabla_{\gamma}\rho\right) = p\left(\frac{Q-p\theta}{p}\right)^{p}\frac{|x|^{\gamma(\alpha+p)}}{\rho^{\gamma(\alpha+p)+p\theta}}
\\&=p\left(\frac{Q-p\theta}{p}\right)^{p}\frac{|\nabla_{\gamma}\rho|^{\alpha+p}}{\rho^{p\theta}}=pw + 0.
\end{align*}
The condition (\ref{main condt}) is satisfied with $\phi=0$, which implies the desired results (\ref{cor4.5 eq0}) and (\ref{cor4.5 eq}). The function $h=\rho^{\frac{Q-p\theta}{p}}$ satisfies the H{\"o}lder equality condition,
\begin{align*}
\left(\frac{p}{Q-p\theta}\right)^{p}\frac{|\nabla_\gamma \rho|^{\alpha}}{\rho^{p(\theta - 1)}}\frac{|\nabla_{\gamma}\rho\cdot\nabla_{\gamma}h|^p}{|\nabla_{\gamma}\rho|^p} = \frac{\left| \nabla_\gamma \rho \right|^{\alpha + p}}{\rho^{p\theta}}|h|^{p},    
\end{align*}
implying the sharpness in the inequality (\ref{cor4.5 eq0}). \QEDA

\vspace{2mm}

Now consider the following pair: 

\begin{align*}
v=\left(\log \frac{R}{\rho}\right)^{\alpha+p},\quad w=\left(\frac{|\alpha+1|}{p}\right)^{p}\left(\log \frac{R}{\rho}\right)^{\alpha}\frac{|x|^{\gamma p}}{\rho^{\gamma p+p}}.
\end{align*}
Then, we are able to obtain the following refinements of the power logarithmic $L^p$-Hardy type inequality \cite[Formula (3.8) of Theorem 3.1]{D'A05}:
\begin{cor}\label{cor4}
Let $1<p<\infty$, $\alpha\in\mathbb{R}$ and $R>0$ be such that $\alpha+1<0$.  
\begin{enumerate}
    \item Then, for all complex-valued $f$ $\in$  $C^{\infty}_{0}(B_{R}\backslash\{(0,0)\})$, there holds
    \begin{align}\label{cor4 eq0}
    \int_{B_{R}}\left(\log\frac{R}{\rho}\right)^{\alpha+p}\frac{|\nabla_{\gamma}\rho\cdot\nabla_{\gamma}f|^p}{|\nabla_{\gamma}\rho|^p}dz\geq\left(\frac{|\alpha+1|}{p}\right)^{p}\int_{B_{R}}\left(\log \frac{R}{\rho}\right)^{\alpha}\frac{|x|^{\gamma p}}{\rho^{\gamma p+p}}|f|^{p}dz,
    \end{align}
    where the constant $\left(\frac{|\alpha+1|}{p}\right)^{p}$ is sharp.
    \item Furthermore, for all complex-valued $f$ $\in$  $C^{\infty}_{0}(B_{R}\backslash\{(0,0)\})$, we also have the identity
\begin{multline}\label{cor4 eq}
\int_{B_{R}}\left(\log\frac{R}{\rho}\right)^{\alpha+p}\frac{|\nabla_{\gamma}\rho\cdot\nabla_{\gamma}f|^p}{|\nabla_{\gamma}\rho|^p}dz=\left(\frac{|\alpha+1|}{p}\right)^{p}\int_{B_{R}}\left(\log \frac{R}{\rho}\right)^{\alpha}\frac{|x|^{\gamma p}}{\rho^{\gamma p+p}}|f|^{p}dz\\+\int_{B_{R}}C_{p}\Biggl(\left(\log\frac{R}{\rho}\right)^{\frac{\alpha+p}{p}}\frac{\nabla_{\gamma}\rho \cdot \nabla_{\gamma}f}{|\nabla_{\gamma}\rho|},\left(\log\frac{R}{\rho}\right)^{\frac{\alpha+p}{p}}\frac{\nabla_{\gamma}\rho \cdot \nabla_{\gamma}f}{|\nabla_{\gamma}\rho|} \\+\left(\frac{|\alpha+1|}{p}\right)\left(\log\frac{R}{\rho}\right)^{\frac{\alpha}{p}}\frac{|x|^{\gamma }}{\rho^{\gamma+1}}f\Biggr)dz
\\+\left(\frac{|\alpha+1|}{p}\right)^{p-1}(Q-p)\int_{B_{R}}\frac{|x|^{\gamma p}}{\rho^{\gamma p + p}}\left(\log \frac{R}{\rho}\right)^{\alpha+1}|f|^{p}dz
\end{multline}
with functional $C_{p}(\cdot,\cdot)$ given in Theorem \ref{thm1}.
\end{enumerate}
\end{cor}
\begin{rem}
If we choose $\gamma=0$, $\alpha=-n$, and $p=n$ in (\ref{cor4 eq}), then for $n\geq2$, this leads to a scale-invariant critical form of the Hardy inequality, as established by Ioku, Ishiwata, and Ozawa \cite{ioku2015scale, ioku2016sharp}:
\begin{multline*}
\int_{B_{R}}\frac{|z\cdot\nabla f|^{n}}{|z|^{n}}dz=\left(\frac{n-1}{n}\right)^{n}\int_{B_R}\frac{|f|^{n}}{|z|^{n}\left(\log \frac{R}{|z|}\right)^n}dz\\+\int_{B_{R}}C_{n}\left(\frac{z\cdot\nabla f}{|z|},\frac{z\cdot\nabla f}{|z|}+\left(\frac{n-1}{n}\right)\frac{f}{|z|\left(\log \frac{R}{|z|}\right)}\right)dz
\end{multline*}
for all complex-valued $f$ $\in$  $C^{\infty}_{0}(B_{R}\backslash\{(0,0)\})$.
\end{rem}
\begin{flushleft}
\textit{Proof of Corollary \ref{cor4}:} Calculating divergence:
\end{flushleft}
\begin{align}
&\nabla_{\gamma}\cdot \left(\left(\left(\frac{|\alpha+1|}{p}\right)^{p}\left(\log \frac{R}{\rho}\right)^{\alpha}\frac{|x|^{\gamma p}}{\rho^{\gamma p+p}}\right)^{\frac{p-1}{p}}\left(\left(\log \frac{R}{\rho}\right)^{\alpha+p}\right)^{\frac{1}{p}}\frac{\nabla_{\gamma}\rho}{|\nabla_{\gamma}\rho|}\right) \nonumber
\\&=\nabla_{\gamma}\cdot\left(\Phi \frac{\nabla_{\gamma}\rho}{|\nabla_{\gamma}\rho|}\right) \nonumber
\\&=\nabla_{\gamma}\Phi \cdot \frac{\nabla_{\gamma}\rho}{|\nabla_{\gamma}\rho|}+\Phi \nabla_{\gamma}\cdot \left(\frac{\nabla_{\gamma}\rho}{|\nabla_{\gamma}\rho|}\right) \nonumber
\\&=\frac{1}{|\nabla_{\gamma}\rho|}\nabla_{\gamma}\Phi \cdot \nabla_{\gamma}\rho + \Phi\nabla_{\gamma}\cdot \left(\frac{\nabla_{\gamma}\rho}{|\nabla_{\gamma}\rho|}\right) \nonumber
\\&=T_{1}+T_{2}, \label{cor4 t1+t2}
\end{align}
where
\begin{align}
&T_1=\frac{1}{|\nabla_{\gamma}\rho|}\nabla_{\gamma}\Phi \cdot \nabla_{\gamma}\rho, \nonumber
\\&T_2=\Phi\nabla_{\gamma}\cdot \left(\frac{\nabla_{\gamma}\rho}{|\nabla_{\gamma}\rho|}\right) \label{cor4 t2}
\end{align}
with
\begin{align*}
\Phi=\left(\frac{|\alpha+1|}{p}\right)^{p-1}\left(\log \frac{R}{\rho}\right)^{\alpha+1}\frac{|x|^{\gamma(p-1)}}{\rho^{(\gamma+1)(p-1)}}.
\end{align*}
By using formula (\ref{formula for term2}) in (\ref{cor4 t2}), we get
\begin{align}
&T_{2}=\Phi\cdot\left[(Q-1)\frac{|x|^{\gamma}}{\rho^{\gamma+1}}\right] \nonumber
\\&=\left(\frac{|\alpha+1|}{p}\right)^{p-1}\left(\log \frac{R}{\rho}\right)^{\alpha+1}\frac{|x|^{\gamma(p-1)}}{\rho^{(\gamma+1)(p-1)}}\cdot(Q-1)\frac{|x|^{\gamma}}{\rho^{\gamma+1}}
\nonumber
\\&=\left(\frac{|\alpha+1|}{p}\right)^{p-1}(Q-1)\left(\log \frac{R}{\rho}\right)^{\alpha+1}\frac{|x|^{\gamma p}}{\rho^{(\gamma+1)p}}. \label{cor4 t2_1}
\end{align}
Now we compute $T_1$:
\begin{align}
&T_{1}=\frac{1}{|\nabla_{\gamma}\rho|}(\nabla_{\gamma}\Phi)\cdot\nabla_{\gamma}\rho \nonumber
\\&=\frac{1}{|\nabla_{\gamma}\rho|}\left(\frac{|\alpha+1|}{p}\right)^{p}\nabla_{\gamma}\left(\left(\log\frac{R}{\rho}\right)^{\alpha+1}|x|^{\gamma(p-1)}\rho^{-(\gamma+1)(p-1)}\right)\cdot\nabla_{\gamma}\rho \nonumber
\\&=-\frac{\rho^{\gamma}}{|x|^{\gamma}}\left(\frac{|\alpha+1|}{p}\right)^{p}\frac{|x|^{\gamma(p+1)}}{\rho^{p(\gamma+1)+\gamma}}\left(\log \frac{R}{\rho}\right)^{\alpha}\left[(\alpha+1)+(p-1)\log \frac{R}{\rho}\right] \nonumber
\\&=-\left(\frac{|\alpha+1|}{p}\right)^{p-1}\frac{|x|^{\gamma p}}{\rho^{p(\gamma+1)}}\left(\log \frac{R}{\rho}\right)^{\alpha}\left[(\alpha+1)+(p-1)\log \frac{R}{\rho}\right]. \label{cor4 t1_1}
\end{align}
Putting (\ref{cor4 t2_1}) and (\ref{cor4 t1_1}) to (\ref{cor4 t1+t2}), one obtains
\allowdisplaybreaks
\begin{align*}
&T_{1}+T_{2}=-\left(\frac{|\alpha+1|}{p}\right)^{p-1}\frac{|x|^{\gamma p}}{\rho^{p(\gamma+1)}}\left(\log \frac{R}{\rho}\right)^{\alpha}\left[(\alpha+1)+(p-1)\log \frac{R}{\rho}\right] \nonumber
\\&+\left(\frac{|\alpha+1|}{p}\right)^{p-1}(Q-1)\left(\log \frac{R}{\rho}\right)^{\alpha+1}\frac{|x|^{\gamma p}}{\rho^{(\gamma+1)p}} \nonumber
\\&=\left(\frac{|\alpha+1|}{p}\right)^{p-1}\frac{|x|^{\gamma p}}{\rho^{(\gamma+1)p}}\left(\log \frac{R}{\rho}\right)^{\alpha}\left[(Q-1)\log \frac{R}{\rho}-(p-1)\log \frac{R}{\rho}-(\alpha+1)\right] \nonumber
\\&=\left(\frac{|\alpha+1|}{p}\right)^{p-1}\frac{|x|^{\gamma p}}{\rho^{(\gamma+1)p}}\left(\log \frac{R}{\rho}\right)^{\alpha}\left[(Q-p)\log \frac{R}{\rho}-(\alpha+1)\right] \nonumber
\\&=\left(\frac{|\alpha+1|}{p}\right)^{p-1}\frac{|x|^{\gamma p}}{\rho^{(\gamma+1)p}}\left(\log \frac{R}{\rho}\right)^{\alpha}\left[(Q-p)\log \frac{R}{\rho}+|\alpha+1|\right] \nonumber
\\&=p\left(\frac{|\alpha+1|}{p}\right)^{p}\frac{|x|^{\gamma p}}{\rho^{\gamma p+p}}\left(\log \frac{R}{\rho}\right)^{\alpha} \nonumber
\\&+\left(\frac{|\alpha+1|}{p}\right)^{p-1}(Q-p)\frac{|x|^{\gamma p}}{\rho^{\gamma p+p}}\left(\log \frac{R}{\rho}\right)^{\alpha+1} \nonumber
\\&=pw+\left(\frac{|\alpha+1|}{p}\right)^{p-1}(Q-p)\frac{|x|^{\gamma p}}{\rho^{\gamma p+p}}\left(\log \frac{R}{\rho}\right)^{\alpha+1}. 
\end{align*}
Therefore, $\phi=\left(\frac{|\alpha+1|}{p}\right)^{p-1}(Q-p)\frac{|x|^{\gamma p}}{\rho^{\gamma p+p}}\left(\log \frac{R}{\rho}\right)^{\alpha+1}\geq0$. The condition (\ref{main condt}) is satisfied and the results (\ref{cor4 eq0}) and (\ref{cor4 eq}) are obtained. Since the function $h=\left(\log \frac{R}{\rho}\right)^{-\frac{|\alpha+1|}{p}}$ satisfies the H{\"o}lder equality condition,
\begin{align*}
\left(\frac{p}{|\alpha+1|}\right)^{p}\left(\log\frac{R}{\rho}\right)^{\alpha+p}\frac{|\nabla_{\gamma}\rho\cdot\nabla_{\gamma}h|^p}{|\nabla_{\gamma}\rho|^p} = \left(\log \frac{R}{\rho}\right)^{\alpha}\frac{|x|^{\gamma p}}{\rho^{\gamma p+p}}|h|^{p},
\end{align*}
the constant in the inequality (\ref{cor4 eq0}) is sharp. \QEDA

\subsection{Caffarelli-Kohn-Nirenberg and Heisenberg-Pauli-Weyl type inequalities}\label{subsection2}

\hspace{1mm}

\begin{flushleft}
In this section, we prove general weighted CKN type inequalities with explicit constants and remainder terms. As a result, we also obtain HPW type inequalities.    
\end{flushleft}

\begin{cor}\label{cor5 ckn}
Let $\Omega$, $\Sigma$, $v$ and $w$ be from Theorem \ref{thm1}. Let $1<p,q<\infty$, $0<r<\infty$ with $p+q\geq r$, $\delta \in [0,1] \cap \left[ \frac{r - q}{r}, \frac{p}{r} \right]$ and $b, c \in \mathbb{R}$. Assume that $\frac{\delta r}{p} + \frac{(1-\delta) r}{q} = 1$ and $c = \frac{\delta}{p} + b(1 - \delta)$. Then we have the following Caffarelli-Kohn-Nirenberg type inequalities for any complex-valued functions $f\in C^{\infty}_{0}(\Omega\backslash\Sigma)$:
\begin{multline}\label{CKN res}
\Biggl[\norm{v^{\frac{1}{p}}\frac{\nabla_{\gamma}\rho\cdot\nabla_{\gamma}f}{|\nabla_{\gamma}\rho|}}^{p}_{L^{p}(\Omega)}\\-\int_{\Omega}C_{p}\left(v^{\frac{1}{p}}\frac{\nabla_{\gamma}\rho \cdot \nabla_{\gamma}f}{|\nabla_{\gamma}\rho|},v^{\frac{1}{p}}\frac{\nabla_{\gamma}\rho \cdot \nabla_{\gamma}f}{|\nabla_{\gamma}\rho|} + w^{\frac{1}{p}}f\right)dz\Biggr]^{\frac{\delta}{p}}\norm{w^{b}f}^{1-\delta}_{L^{q}(\Omega)}\geq\norm{w^{c}f}_{L^{r}(\Omega)}.
\end{multline}  
\end{cor}
\begin{flushleft}
\textit{Proof of Corollary \ref{cor5 ckn}:} Case $\delta=0$. This gives us $q=r$ and $b=c$ by $\frac{\delta r}{p}+\frac{(1-\delta)r}{q}=1$ and $c=\frac{\delta}{p}+b(1-\delta)$, respectively. Then, the inequality (\ref{CKN res}) reduces to the trivial estimate
\end{flushleft}
\begin{align*}
\norm{w^{b}f}_{L^{q}(\Omega)}\leq\norm{w^{{b}}f}_{L^{q}(\Omega)}.
\end{align*}
Case $\delta=1$. In this case, we get $p=r$ and $c=\frac{1}{p}$ with
\begin{multline*}
\norm{w^{\frac{1}{p}}f}^{p}_{L^{p}(\Omega)}\leq\norm{v^{\frac{1}{p}}\frac{\nabla_{\gamma}\rho\cdot\nabla_{\gamma}f}{|\nabla_{\gamma}\rho|}}^{p}_{L^{p}(\Omega)}\\-\int_{\Omega}C_{p}\left(v^{\frac{1}{p}}\frac{\nabla_{\gamma}\rho \cdot \nabla_{\gamma}f}{|\nabla_{\gamma}\rho|},v^{\frac{1}{p}}\frac{\nabla_{\gamma}\rho \cdot \nabla_{\gamma}f}{|\nabla_{\gamma}\rho|} + w^{\frac{1}{p}}f\right)dz,
\end{multline*}
which is (\ref{main res}) by Theorem \ref{thm1}.
\begin{flushleft}
\vspace{1mm}
Case $\delta\in (0,1)\cap\left[\frac{r-q}{q},\frac{p}{r}\right]$. Taking into consideration $c=\frac{\delta}{p}+b(1-\delta)$, we can directly obtain 
\end{flushleft}
\begin{align*}
\norm{w^{c}f}_{L^{r}(\Omega)}=\left(\int_{\Omega}w^{cr}|f|^rdz\right)^{\frac{1}{r}}=\left(\int_{\Omega}\frac{|f|^{\delta r}}{w^{-\frac{\delta r}{p}}}\cdot\frac{|f|^{(1-\delta) r}}{w^{-b(1-\delta)r}}dz\right)^{\frac{1}{r}}.
\end{align*}
Since $\delta\in (0,1)\cap\left[\frac{r-q}{q},\frac{p}{r}\right]$ and $p+q\geq r$, then using the H{\"o}lder inequality for $\frac{\delta r}{p} + \frac{(1-\delta) r}{q} = 1$, we derive
\begin{align}
\norm{w^{c}f}_{L^{r}(\Omega)}&\leq\left(\int_{\Omega}w|f|^{p}\right)^\frac{\delta}{p}\left(\int_{\Omega}w^{bq}|f|^{q}\right)^{\frac{1-\delta}{q}} \nonumber
\\&=\norm{w^{\frac{1}{p}}f}_{L^{p}(\Omega)}^{\delta}\norm{w^{b}f}^{1-\delta}_{L^{q}(\Omega)} \nonumber
\\&\leq\Biggl[\norm{v^{\frac{1}{p}}\frac{\nabla_{\gamma}\rho\cdot\nabla_{\gamma}f}{|\nabla_{\gamma}\rho|}}^{p}_{L^{p}(\Omega)} \nonumber
\\&-\int_{\Omega}C_{p}\left(v^{\frac{1}{p}}\frac{\nabla_{\gamma}\rho \cdot \nabla_{\gamma}f}{|\nabla_{\gamma}\rho|},v^{\frac{1}{p}}\frac{\nabla_{\gamma}\rho \cdot \nabla_{\gamma}f}{|\nabla_{\gamma}\rho|} + w^{\frac{1}{p}}f\right)dz\Biggr]^{\frac{\delta}{p}}\norm{w^{b}f}^{1-\delta}_{L^{q}(\Omega)}.
\label{ckn last eq}
\end{align}
Rewriting (\ref{ckn last eq}), we derive (\ref{CKN res}). \QEDA

\vspace{3mm}
By the same technique, we can apply Corollaries \ref{cor1}, \ref{cor2}, \ref{cor4.5} and \ref{cor4} to obtain the corresponding CKN type inequalities. 

\begin{cor}
Let $1<p,q<\infty$, $0<r<\infty$ with $p+q\geq r$, $\delta \in [0,1] \cap \left[ \frac{r - q}{r}, \frac{p}{r} \right]$ and $b, c \in \mathbb{R}$. Assume that $\frac{\delta r}{p} + \frac{(1-\delta) r}{q} = 1$ and $c = \frac{\delta}{p} + b(1 - \delta)$. Then,
\begin{enumerate}
    \item for any complex-valued $f\in C^{\infty}_{0}(B_{R}\backslash\{(0,0)\})$, we have
    \begin{multline}\label{hpw cor1}
\Biggl[\norm{\frac{\nabla_{\gamma}\rho\cdot\nabla_{\gamma}f}{|\nabla_{\gamma}\rho|}}^{p}_{L^{p}(B_{R})}
-\int_{B_R}C_{p}\biggl(\frac{\nabla_{\gamma}\rho \cdot \nabla_{\gamma}f}{|\nabla_{\gamma}\rho|},\frac{\nabla_{\gamma}\rho \cdot \nabla_{\gamma}f}{|\nabla_{\gamma}\rho|}\\ + \left(\left(\frac{p-1}{p}\right)^p\frac{|x|^{\gamma p}}{(R-\rho)^p\rho^{\gamma p}}\right)^{\frac{1}{p}}f\biggr)dz\Biggr]^{\frac{\delta}{p}}\norm{\frac{|x|^{\gamma p b}}{(R-\rho)^{pb}\rho^{\gamma p b}}f}^{1-\delta}_{L^{q}(B_{R})}\\\geq\left(\frac{p-1}{p}\right)^{\delta}\norm{\frac{|x|^{\gamma pc}}{(R-\rho)^{pc}\rho^{\gamma pc}}f}_{L^{r}(B_{R})};
\end{multline}
\item for any complex-valued $f\in C^{\infty}_{0}(\mathbb{R}^{m}\times\mathbb{R}^{k}\backslash\{(0,0)\})$, we have
\begin{multline}\label{hpw cor2}
\Biggl[\norm{|x|^{\frac{\beta-\gamma p}{p}}\rho^{\frac{p(1+\gamma)-\alpha}{p}}\frac{\nabla_{\gamma}\rho\cdot\nabla_{\gamma}f}{|\nabla_{\gamma}\rho|}}^{p}_{L^{p}(\mathbb{R}^{n})}\\-\int_{\mathbb{R}^n}C_{p}\biggl(|x|^{\frac{\beta-\gamma p}{p}}\rho^{\frac{p(1+\gamma)-\alpha}{p}}\frac{\nabla_{\gamma}\rho \cdot \nabla_{\gamma}f}{|\nabla_{\gamma}\rho|},|x|^{\frac{\beta-\gamma p}{p}}\rho^{\frac{p(1+\gamma)-\alpha}{p}}\frac{\nabla_{\gamma}\rho \cdot \nabla_{\gamma}f}{|\nabla_{\gamma}\rho|} \\+ \left(\frac{Q+\beta-\alpha}{p}\right)\frac{|x|^{\frac{\beta}{p}}}{\rho^{\frac{\alpha}{p}}}f\biggr)dz\Biggr]^{\frac{\delta}{p}}\norm{\frac{|x|^{\beta b}}{\rho^{\alpha b}}f}^{1-\delta}_{L^{q}(\mathbb{R}^n)}\\\geq\left(\frac{Q+\beta-\alpha}{p}\right)^{\delta}\norm{\frac{|x|^{\beta c}}{\rho^{\alpha c}}f}_{L^{r}(\mathbb{R}^{n})};
\end{multline}
\item for any complex-valued $f$ $\in$  $C^{\infty}_{0}(B_{R}\backslash\{(0,0)\})$, we have
\begin{multline*}
\Biggl[\norm{\frac{|\nabla_\gamma \rho|^{\frac{\alpha}{p}}}{\rho^{(\theta - 1)}}\frac{\nabla_{\gamma}\rho \cdot \nabla_{\gamma}f}{|\nabla_{\gamma}\rho|}}^{p}_{L^{p}(B_{R})}-\int_{B_{R}}C_{p}\biggl(\frac{|\nabla_\gamma \rho|^{\frac{\alpha}{p}}}{\rho^{(\theta - 1)}}\frac{\nabla_{\gamma}\rho \cdot \nabla_{\gamma}f}{|\nabla_{\gamma}\rho|},\frac{|\nabla_\gamma \rho|^{\frac{\alpha}{p}}}{\rho^{(\theta - 1)}}\frac{\nabla_{\gamma}\rho \cdot \nabla_{\gamma}f}{|\nabla_{\gamma}\rho|} \\+ \left(\frac{Q - p\theta}{p}\right)\frac{\left| \nabla_\gamma \rho \right|^{\frac{\alpha + p}{p}}}{\rho^{\theta}}f\biggr)dz\Biggr]^{\frac{\delta}{p}}\norm{\frac{\left| \nabla_\gamma \rho \right|^{\alpha b + pb}}{\rho^{p\theta b}}f}^{1-\delta}_{L^{q}(B_{R})}\\\geq\left(\frac{Q-p\theta}{p}\right)^{\delta}\norm{\frac{\left| \nabla_\gamma \rho \right|^{\alpha c + pc}}{\rho^{p\theta c}}f}_{L^{r}(B_{R})};
\end{multline*}
\item for any complex-valued $f$ $\in$  $C^{\infty}_{0}(B_{R}\backslash\{(0,0)\})$, we have
\begin{multline}\label{hpw cor3}
\Biggl[\norm{\left(\log\frac{R}{\rho}\right)^{\frac{\alpha+p}{p}}\frac{\nabla_{\gamma}\rho\cdot\nabla_{\gamma}f}{|\nabla_{\gamma}\rho|}}^{p}_{L^{p}(B_{R})}\\-\int_{B_{R}}C_{p}\Biggl(\left(\log\frac{R}{\rho}\right)^{\frac{\alpha+p}{p}}\frac{\nabla_{\gamma}\rho \cdot \nabla_{\gamma}f}{|\nabla_{\gamma}\rho|},\left(\log\frac{R}{\rho}\right)^{\frac{\alpha+p}{p}}\frac{\nabla_{\gamma}\rho \cdot \nabla_{\gamma}f}{|\nabla_{\gamma}\rho|} \\+ \left(\frac{|\alpha+1|}{p}\right)\left(\log\frac{R}{\rho}\right)^{\frac{\alpha}{p}}\frac{|x|^{\gamma }}{\rho^{\gamma
+1}}f\Biggr)dz\Biggr]^{\frac{\delta}{p}}\norm{\left(\log \frac{R}{\rho}\right)^{\alpha b}\frac{|x|^{\gamma p b}}{\rho^{\gamma pb+pb}}f}^{1-\delta}_{L^{q}(B_{R})}\\\geq\left(\frac{|\alpha+1|}{p}\right)^{\delta}\norm{\left(\log \frac{R}{\rho}\right)^{\alpha c}\frac{|x|^{\gamma pc}}{\rho^{\gamma pc+pc}}f}_{L^{r}(B_{R})}.
\end{multline}
\end{enumerate}
\end{cor}

In the special case, we also obtain some HPW type inequalities.

\begin{cor}\label{cor9}
Let $p'=\frac{p}{p-1}$. Then,
\begin{enumerate}
    \item for all complex-valued $f\in C^{\infty}_{0}(B_{R}\backslash\{(0,0)\})$, we have
    \begin{align}\label{hpw cor1 p'}
    \left(\int_{B_{R}}\frac{|\nabla_{\gamma}\rho\cdot\nabla_{\gamma}f|^{p}}{|\nabla_{\gamma}\rho|^{p}}dz\right)^{\frac{1}{p}}\left(\int_{B_{R}}\frac{(R-\rho)^{\frac{p p'}{2}}\rho^{\frac{\gamma p p'}{2}}}{|x|^{\frac{\gamma pp'}{2}}}|f|^{p'}dz\right)^{\frac{1}{p'}}\geq\left(\frac{p-1}{p}\right)\left(\int_{B_{R}}|f|^{2}dz\right);
    \end{align}
    \item for all complex-valued $f\in C^{\infty}_{0}(\mathbb{R}^m\times\mathbb{R}^{n}\backslash\{(0,0)\})$, we have
    \begin{align}\label{hpw dam p'}
    \left(\int_{\mathbb{R}^n}\frac{|\nabla_{\gamma}\rho\cdot\nabla_{\gamma}f|^{p}}{|\nabla_{\gamma}\rho|^{p}}dz\right)^{\frac{1}{p}}\left(\int_{\mathbb{R}^n}\frac{\rho^{\frac{\gamma p p'+pp'}{2}}}{|x|^{\frac{\gamma p p'}{2}}}|f|^{p'}dz\right)^{\frac{1}{p'}}\geq\left(\frac{Q-p}{p}\right)\left(\int_{\mathbb{R}^n}|f|^{2}dz\right);
    \end{align}
    \item for all complex-valued $f\in C^{\infty}_{0}(B_{R}\backslash\{(0,0)\})$, we have
    \begin{multline}\label{hpw log p'}
    \left(\int_{B_{R}}\left(\log \frac{R}{\rho}\right)^{2p}\frac{|\nabla_{\gamma}\rho\cdot\nabla_{\gamma}f|^{p}}{|\nabla_{\gamma}\rho|^{p}}dz\right)^{\frac{1}{p}}\left(\int_{B_{R}}\frac{\rho^{\frac{\gamma p p'+ pp'}{2}}}{|x|^{\frac{\gamma p p'}{2}}\left(\log\frac{R}{\rho}\right)^{\frac{pp'}{2}}}|f|^{p'}dz\right)^{\frac{1}{p'}}\\\geq\left(\frac{|p+1|}{p}\right)\left(\int_{B_{R}}|f|^{2}dz\right).
    \end{multline}
\end{enumerate}
\end{cor}

\begin{rem}
In the case when $p=p'=2$ in (\ref{hpw dam p'}), we get
\begin{align*}
\left(\int_{\mathbb{R}^n}\frac{|\nabla_{\gamma}\rho\cdot\nabla_{\gamma}f|^{2}}{|\nabla_{\gamma}\rho|^{2}}dz\right)\left(\int_{\mathbb{R}^n}\frac{\rho^{2\gamma+2}}{|x|^{2\gamma}}|f|^{2}dz\right)\geq\left(\frac{Q-2}{2}\right)^{2}\left(\int_{\mathbb{R}^n}|f|^{2}dz\right)^{2}.
\end{align*}
By choosing $\gamma=0$ and using the Cauchy-Schwarz inequality, we arrive at the following uncertainty principle with an explicit constant:
\begin{align*}
\left( \int_{\mathbb{R}^n} |\nabla f|^2 \, dz \right)
\left( \int_{\mathbb{R}^n} |z|^2 |f|^2 \, dz \right)
\geq \frac{(n-2)^2}{4} \left( \int_{\mathbb{R}^n} |f|^2 \, dz \right)^2.
\end{align*}
\end{rem}

\begin{flushleft}
\textit{Proof of Corollary \ref{cor9}:} Let $c=0$, $r=2$, $b=-\frac{1}{2}$, $q=p'=\frac{p}{p-1}$, $\delta=\frac{1}{2}$. Then, by substituting these values into (\ref{hpw cor1}), we get (\ref{hpw cor1 p'}). Furthermore, taking  $\beta=\gamma p$ and $\alpha=p(1+\gamma)$ in (\ref{hpw cor2}), we obtain (\ref{hpw dam p'}). Finally, choosing $\alpha=p$ in (\ref{hpw cor3}), we derive (\ref{hpw log p'}). \QEDA
\end{flushleft}

{\bf Data availability} Data sharing not applicable to this article as no datasets were generated or analyzed during
the current study

{\bf Declarations}

{\bf Conflict of interest} The authors declared that they have no conflict of interest to this work.

\bibliographystyle{alpha}
\bibliography{citation}

\begin{thebibliography}{DEFT15}

\bibitem[BEL15]{BEL15}
A.~A. Balinsky, W.~D. Evans, and R.~T. Lewis.
\newblock {\em The analysis and geometry of {H}ardy's inequality}, volume~1.
\newblock Springer, 2015.

\bibitem[CFLL24]{cazacu2024caffarelli}
C.~Cazacu, J.~Flynn, N.~Lam, and G.~Lu.
\newblock Caffarelli-{K}ohn-{N}irenberg identities, inequalities and their stabilities.
\newblock {\em J. Math. Pures Appl.}, 182:253--284, 2024.

\bibitem[CKLL24]{cazacu2024hardy}
C.~Cazacu, D.~Krej{\v{c}}i{\v{r}}{\'\i}k, N.~Lam, and A.~Laptev.
\newblock Hardy inequalities for magnetic $p$-{L}aplacians.
\newblock {\em Nonlinearity}, 37(3):035004, 2024.

\bibitem[CKN84]{caffarelli1984first}
L.~Caffarelli, R.~Kohn, and L.~Nirenberg.
\newblock First order interpolation inequalities with weights.
\newblock {\em Compos. Math.}, 53(3):259--275, 1984.

\bibitem[CT24]{CT24}
X.~P. Chen and C.~L. Tang.
\newblock Remainder terms of {$L^p$}-{H}ardy inequalities with magnetic fields: the case $1<p<2$.
\newblock {\em 10.48550/arXiv.2408.17249}, 2024.

\bibitem[D'A04]{D'A04}
L.~D'Ambrosio.
\newblock Hardy inequalities related to {G}rushin type operators.
\newblock {\em Proc. Amer. Math. Soc.}, 132(3):725--734, 2004.

\bibitem[D'A05]{D'A05}
L.~D'Ambrosio.
\newblock Hardy-type inequalities related to degenerate elliptic differential operators.
\newblock {\em Ann. Sc. Norm. Super. Pisa. Cl. Sci.}, 4(3):451--486, 2005.

\bibitem[D'A24]{d2024weighted}
L.~D'Arca.
\newblock Weighted {P}oincar\'e inequality and {H}ardy improvements related to some degenerate elliptic differential operators.
\newblock {\em arXiv preprint arXiv:2407.10840}, 2024.

\bibitem[DEFT15]{dolbeault2015rigidity}
J.~Dolbeault, M.~J. Esteban, S.~Filippas, and A.~Tertikas.
\newblock Rigidity results with applications to best constants and symmetry of {C}affarelli-{K}ohn-{N}irenberg and logarithmic {H}ardy inequalities.
\newblock {\em Calc. Var. Partial Differential Equations}, 54:2465--2481, 2015.

\bibitem[DEL16]{dolbeault2016rigidity}
J.~Dolbeault, M.~J. Esteban, and M.~Loss.
\newblock Rigidity versus symmetry breaking via nonlinear flows on cylinders and {E}uclidean spaces.
\newblock {\em Invent. Math.}, 206:397--440, 2016.

\bibitem[DGN10]{DGN10}
J.~Dou, Q.~Guo, and P.~Niu.
\newblock Hardy inequalities with remainder terms for the generalized {B}aouendi-{G}rushin vector fields.
\newblock {\em Math. Inequal. Appl.}, 13(3):555--570, 2010.

\bibitem[Don18]{dong2018existence}
M.~Dong.
\newblock Existence of extremal functions for higher-order {C}affarelli-{K}ohn-{N}irenberg inequalities.
\newblock {\em Adv. Nonlinear Stud.}, 18(3):543--553, 2018.

\bibitem[Gar93]{Gar93}
N.~Garofalo.
\newblock Unique continuation for a class of elliptic operators which degenerate on a manifold of arbitrary codimension.
\newblock {\em J. Differential Equations}, 104(1):117--146, 1993.

\bibitem[GJR24]{GJR24}
D.~Ganguly, K.~Jotsaroop, and P.~Roychowdhury.
\newblock {H}ardy and {R}ellich identities and inequalities for {B}aouendi-{G}rushin operators via spherical vector fields.
\newblock {\em arXiv preprint arXiv:2404.05510}, 2024.

\bibitem[GL90]{GL90}
N.~Garofalo and E.~Lanconelli.
\newblock Frequency functions on the {H}eisenberg group, the uncertainty principle and unique continuation.
\newblock {\em Ann. Inst. Fourier (Grenoble)}, 40(2):313--356, 1990.

\bibitem[HZ11]{han2011class}
Y.~Z. Han and Q.~Zhao.
\newblock A class of {C}affarelli-{K}ohn-{N}irenberg type inequalities for generalized {B}aouendi-{G}rushin vector fields.
\newblock {\em Acta Math. Sci. Ser. A Chin. Ed}, 31(5):1181--1189, 2011.

\bibitem[II15]{ioku2015scale}
N.~Ioku and M.~Ishiwata.
\newblock A scale invariant form of a critical {H}ardy inequality.
\newblock {\em Int. Math. Res. Not. IMRN}, 2015(18):8830--8846, 2015.

\bibitem[IIO16]{ioku2016sharp}
N.~Ioku, M.~Ishiwata, and T.~Ozawa.
\newblock Sharp remainder of a critical {H}ardy inequality.
\newblock {\em Arch. Math. (Basel)}, 106(1), 2016.

\bibitem[KMP06]{KMP06}
A.~Kufner, L.~Maligranda, and L.~E. Persson.
\newblock The prehistory of the {H}ardy inequality.
\newblock {\em Amer. Math. Monthly}, 113(8):715--732, 2006.

\bibitem[KMP07]{KMP07}
A.~Kufner, L.~Maligranda, and L.~E. Persson.
\newblock {\em The {H}ardy inequality: {A}bout its history and some related results}.
\newblock Vydavatelsk{\`y} servis, 2007.

\bibitem[Kom06]{kombe2006hardy}
I.~Kombe.
\newblock Hardy, {R}ellich and uncertainty principle inequalities on {C}arnot groups.
\newblock {\em arXiv preprint math/0611850}, 2006.

\bibitem[Kom15]{Kom15}
I.~Kombe.
\newblock Hardy and {R}ellich type inequalities with remainders for {B}aouendi-{G}rushin vector fields.
\newblock {\em Houston J. Math}, 41(3), 2015.

\bibitem[KY18]{KY18}
I.~Kombe and A.~Yener.
\newblock General weighted {H}ardy type inequalities related to {B}aouendi-{G}rushin operators.
\newblock {\em Complex Var. Elliptic Equ.}, 63(3):420--436, 2018.

\bibitem[KY24]{kalaman2024cylindrical}
M.~Kalaman and N.~Yessirkegenov.
\newblock Cylindrical {H}ardy, {S}obolev type and {C}affarelli-{K}ohn-{N}irenberg type inequalities and identities.
\newblock {\em arXiv preprint arXiv:2407.08393}, 2024.

\bibitem[LL17]{lam2017sharp}
N.~Lam and G.~Lu.
\newblock Sharp constants and optimizers for a class of {C}affarelli-{K}ohn-{N}irenberg inequalities.
\newblock {\em Adv. Nonlinear Stud.}, 17(3):457--480, 2017.

\bibitem[LRY19]{LRY19}
A.~Laptev, M.~Ruzhansky, and N.~Yessirkegenov.
\newblock Hardy inequalities for {L}andau {H}amiltonian and for {B}aouendi-{G}rushin operator with {A}haronov-{B}ohm type magnetic field. {P}art {I}.
\newblock {\em Math. Scand.}, 125(2):239--269, 2019.

\bibitem[LY23]{li2023anisotropic}
YY. Li and X.~Yan.
\newblock Anisotropic {C}affarelli-{K}ohn-{N}irenberg type inequalities.
\newblock {\em Adv. Math.}, 419:108958, 2023.

\bibitem[NCH04]{NCH04}
P.~Niu, Y.~Chen, and Y.~Han.
\newblock Some {H}ardy-type inequalities for the generalized {B}aouendi-{G}rushin operators.
\newblock {\em Glasg. Math. J.}, 46(3):515--527, 2004.

\bibitem[ORS19]{ozawa2019p}
T.~Ozawa, M.~Ruzhansky, and D.~Suragan.
\newblock ${L}^{p}$-{C}affarelli-{K}ohn--{N}irenberg type inequalities on homogeneous groups.
\newblock {\em Q. J. Math.}, 70(1):305--318, 2019.

\bibitem[RS17a]{ruzhansky2017hardy}
M.~Ruzhansky and D.~Suragan.
\newblock Hardy and {R}ellich inequalities, identities, and sharp remainders on homogeneous groups.
\newblock {\em Adv. Math.}, 317:799--822, 2017.

\bibitem[RS17b]{ruzhansky2017horizontal}
M.~Ruzhansky and D.~Suragan.
\newblock On horizontal {H}ardy, {R}ellich, {C}affarelli-{K}ohn-{N}irenberg and p-sub-{L}aplacian inequalities on stratified groups.
\newblock {\em J. Differential Equations}, 262(3):1799--1821, 2017.

\bibitem[RSY17a]{ruzhansky2017caffarelli}
M.~Ruzhansky, D.~Suragan, and N.~Yessirkegenov.
\newblock Caffarelli-{K}ohn-{N}irenberg and {S}obolev type inequalities on stratified lie groups.
\newblock {\em NoDEA Nonlinear Differential Equations Appl.}, 24(5):56, 2017.

\bibitem[RSY17b]{ruzhansky2017extended}
M.~Ruzhansky, D.~Suragan, and N.~Yessirkegenov.
\newblock Extended {C}affarelli-{K}ohn-{N}irenberg inequalities and superweights for ${L}^{p}$-weighted {H}ardy inequalities.
\newblock {\em C. R. Math. Acad. Sci. Paris}, 355(6):694--698, 2017.

\bibitem[RSY18]{ruzhansky2018extended}
M.~Ruzhansky, D.~Suragan, and N.~Yessirkegenov.
\newblock Extended {C}affarelli-{K}ohn-{N}irenberg inequalities, and remainders, stability, and superweights for ${L}^{p}$-weighted {H}ardy inequalities.
\newblock {\em Trans. Amer. Math. Soc. Ser. B}, 5(2):32--62, 2018.

\bibitem[SJ12]{SJ12}
S.~F. Shen and Y.~Y. Jin.
\newblock Rellich type inequalities related to {G}rushin type operator and {G}reiner type operator.
\newblock {\em Appl. Math. J. Chinese Univ. Ser. A}, 27:353--362, 2012.

\bibitem[SY12]{SY12}
D.~Su and Q.~H. Yang.
\newblock Improved {H}ardy inequalities in the {G}rushin plane.
\newblock {\em J. Math. Anal. Appl.}, 393(2):509--516, 2012.

\bibitem[SY23]{suragan2023sharp}
D.~Suragan and N.~Yessirkegenov.
\newblock Sharp remainder of the {P}oincar{\'e} inequality for {B}aouendi--{G}rushin vector fields.
\newblock {\em Asian-Eur. J. Math.}, 16(03):2350041, 2023.

\bibitem[Wey50]{weyl1950theory}
H.~Weyl.
\newblock {\em The theory of groups and quantum mechanics}.
\newblock Courier Corporation, 1950.

\bibitem[YSK15]{yang2015improved}
Q.~Yang, D.~Su, and Y.~Kong.
\newblock Improved {H}ardy inequalities for {G}rushin operators.
\newblock {\em J. Math. Anal. Appl.}, 424(1):321--343, 2015.

\end{thebibliography}

\end{document}